\documentclass[12pt]{amsart}
\usepackage[a4paper,margin=1in]{geometry}
\usepackage{amsmath,amssymb,amsthm,mathtools}
\usepackage{enumitem}
\usepackage{microtype}
\usepackage{xcolor}
\usepackage[hidelinks]{hyperref}
\newtheorem{theorem}{Theorem}[section]
\newtheorem{proposition}[theorem]{Proposition}
\newtheorem{lemma}[theorem]{Lemma}

\theoremstyle{definition}
\newtheorem{definition}[theorem]{Definition}
\newtheorem{remark}[theorem]{Remark}

\numberwithin{equation}{section}
\newcommand{\R}{\mathbb R}
\newcommand{\cL}{\mathcal L}
\newcommand{\cP}{\mathcal P}

\newcommand{\cQ}{\mathcal Q}
\newcommand{\BMO}{\mathrm{BMO}}
\newcommand{\supp}{\operatorname{supp}}
\newcommand{\dist}{\operatorname{dist}}

\title[Hardy and BMO spaces for the Laplacian with drift]{Hardy and BMO Spaces for the Laplacian with Drift}

\author[T. A. Bui]{The Anh Bui}
\address{School of Mathematical and Physical Sciences, Macquarie University, NSW 2109, Australia}
\email{the.bui@mq.edu.au}

\author{Fu Ken Ly}
\address{School of Mathematics and Statistics, The Learning Hub, The University of Sydney, Sydney, NSW 2006, Australia}
\email{ken.ly@sydney.edu.au}

\subjclass[2020]{Primary 42B30; Secondary 42B20, 42B35}
\keywords{Laplacian with drift, non-doubling measure, Hardy space, BMO space, Riesz transforms, commutators}

\begin{document}
\begin{abstract}
We study Hardy and BMO spaces associated with the
Laplacian with drift. We define an atomic
Hardy space $H^1$ whose atoms satisfy cancellation
against constants and linear functions in directions
perpendicular to the drift. We prove the
John--Nirenberg inequality, Hardy--BMO duality,
equivalence of finite and infinite atomic norms,
and interpolation theorems. We show that the
Riesz transforms and their adjoints are bounded
from $H^1$ to $L^1$ and from $L^\infty$ to BMO.
Moreover, we prove that the cancellation conditions
are optimal, but the Riesz transforms do not
characterize $H^1$. We also characterize the
$L^p$ boundedness of their commutators,
$1<p<\infty$, in terms of local BMO.
\end{abstract}
\maketitle

\tableofcontents
\medskip

\section{Introduction}\label{sec:introduction}

Drift terms appear in many partial differential equations
and differential operators. For drift--diffusion equations,
Caffarelli and Vasseur~\cite{CV} established regularity
results motivated by the quasi-geostrophic equation.
H\"older continuity and potential estimates for related
drift--diffusion equations were studied in~\cite{NST,NNSW}.
For second-order differential operators with drift,
D'Ancona, Maheux and Pierfelice~\cite{DMP} obtained
logarithmic Sobolev inequalities and related heat
kernel estimates. In this paper, for
$v\in\R^n\setminus\{0\}$, we consider the Laplacian
with constant drift
\begin{equation}\label{eq-general L with drift}
	\cL_v=-\Delta-2v\cdot\nabla,
	\qquad
	d\mu_v(x)=e^{2\langle v,x\rangle}\,dx,
\end{equation}
acting on $L^2(\mu_v)$. It is well known that
$\cL_v$ is nonnegative and self-adjoint on
$L^2(\mu_v)$. The Riesz transforms associated
with $\cL_v$ are defined by
\[
R_{v,j}=\partial_j\cL_v^{-1/2},
\qquad 1\le j\le n.
\]

Let $\alpha=|v|$ and choose an orthogonal map $O$
such that $O(v/|v|)=e_1$. By the change of variables
$z=\alpha Ox$, we have
\[
\langle v,x\rangle=z_1,\qquad
d\mu_v(x)=\alpha^{-n}e^{2z_1}\,dz,
\]
and
\[
\cL_v=\alpha^2(-\Delta_z-2\partial_{z_1}).
\]
Thus, by rotation and dilation, we are reduced
to the normalized Laplacian with drift
\begin{equation}\label{eq:L-normalized}
	\cL=-\Delta-2\partial_1,
	\qquad d\mu(x)=e^{2x_1}\,dx.
\end{equation}
Instead of studying $R_{v,j}$, it suffices
to consider
\[
R_j=\partial_j\cL^{-1/2},
\qquad 1\le j\le n.
\]
Indeed, the Riesz transforms associated with
$\cL_v$ are linear combinations of those
associated with $\cL$ after rotation and dilation.

For $f\in C_c^\infty(\R^n)$, we have
\[
\cL f=-e^{-2x_1}
\operatorname{div}(e^{2x_1}\nabla f),
\qquad
\langle\cL f,f\rangle
=\int_{\R^n}|\nabla f|^2\,d\mu,
\]
where $\langle f,g\rangle=\int f\overline g\,d\mu$
is the $L^2(\mu)$ pairing. Since
\[
\partial_j^*=-\partial_j-2\delta_{1j},
\qquad
\cL=\sum_{j=1}^n\partial_j^*\partial_j,
\]
we obtain
\begin{equation}\label{eq:L2-Riesz}
	\sum_{j=1}^n\|R_jf\|_{L^2(\mu)}^2
	=\|f\|_{L^2(\mu)}^2,
	\qquad
	\|R_j\|_{L^2(\mu)\to L^2(\mu)}\le1,
	\qquad
	\|R_j^*\|_{L^2(\mu)\to L^2(\mu)}\le1.
\end{equation}

The Riesz transforms associated with the Laplacian
with drift have been studied in~\cite{LM,LSW,LS-endpoint,Mukherjee},
including $L^p$ estimates and endpoint results.
For Riesz transforms on manifolds, see Coulhon and
Duong~\cite{CD99} and Auscher, Coulhon, Duong
and Hofmann~\cite{ACDH04}.

Harmonic analysis in non-doubling settings has
attracted considerable attention.
{Maximal operators have also been studied in~\cite{AFN,DN}, including the case of rotation-invariant measures. For singular integrals on non-doubling spaces, see~\cite{NTV,YY}.}
Hardy and BMO spaces for Gaussian measures were
studied in~\cite{MM,MMS}, including Hardy--BMO
duality and endpoint estimates for Riesz transforms.
In~\cite{CMM}, the theory of Hardy and BMO spaces
on locally doubling spaces was developed, including
atomic decompositions, Hardy--BMO duality and
interpolation theorems. On Lie groups of exponential
growth, Sj\"ogren and Vallarino~\cite{SV,SV2025}
studied Hardy spaces and endpoint estimates for
Riesz transforms, including the failure of a
Riesz-transform characterization of the atomic
Hardy space. For non-doubling measures satisfying
polynomial growth conditions, Hardy and BMO type
spaces were studied in~\cite{MMNO,Tolsa}.

The main aim of this paper is to establish the
theory of Hardy and BMO spaces associated with
the Laplacian with drift. We also study the
boundedness of Riesz transforms and their
commutators on these spaces. We first consider
the normalized operator \eqref{eq:L-normalized}
and then extend the results to the general
operator \eqref{eq-general L with drift}.

One of the difficulties is that the measure
$\mu$ is not doubling. In fact, $\mu$ is
doubling on balls of bounded radius, but
it grows exponentially in the direction
of the drift. Although the Riesz transforms
are bounded on $L^p(\mu)$ for $1<p<\infty$,
their boundedness from Hardy spaces to
$L^1(\mu)$ requires more cancellation
conditions than in the classical case.

Recall that, for the classical Hardy space
$H^1(\R^n)$, an atom $a$ requires only
the cancellation condition
\[
\int_{\R^n}a(x)\,dx=0.
\]
However, this condition is not sufficient
in our setting. More precisely, we require
our atoms to satisfy
\[
\int a\,d\mu=0,\qquad
\int x_ka(x)\,d\mu(x)=0,
\quad 2\le k\le n.
\]
Thus, instead of requiring cancellation
against constants, we require cancellation
against the space
\[
\cP=\operatorname{span}_{\mathbb C}
\{1,x_2,\ldots,x_n\}.
\]
See Definition~\ref{def:H1} for the precise
definition of our atoms.

The additional cancellation conditions
are needed to control the Riesz transforms
away from the support of an atom.
Using these atoms, we define the Hardy space
$H^1(\mu)$ by atomic decompositions.
We prove that, for $1\le j\le n$,
\begin{equation}\label{eq:intro-H1-endpoint}
	\|R_jf\|_{L^1(\mu)}
	+\|R_j^*f\|_{L^1(\mu)}
	\lesssim\|f\|_{H^1(\mu)}.
\end{equation}
See Theorem~\ref{thm:H1-endpoint}.

More importantly, these cancellation conditions
are optimal. In Theorem~\ref{thm:sharp-moment},
we show that if $\{1,x_2,\ldots,x_n\}$ is replaced
by any proper subset, then $R_1$ is no longer
bounded from the corresponding atomic Hardy
space to $L^1(\mu)$. 

We also establish interpolation theorems
for $H^1(\mu)$. In particular, if a sublinear
operator $T$ is bounded from $H^1(\mu)$ to
$L^1(\mu)$ and on $L^2(\mu)$, then $T$ is
bounded on $L^p(\mu)$ for every $1<p<2$.
See Theorem~\ref{thm:Hardy-interpolation}.

On the other hand, the Riesz transforms
do not characterize $H^1(\mu)$. Indeed,
Theorem~\ref{thm:failure-Riesz-characterization}
shows that there exists $f\notin H^1(\mu)$
such that
\[
f,\ R_jf,\ R_j^*f\in L^1(\mu),
\qquad 1\le j\le n.
\]
This is in contrast to the classical
Riesz-transform characterization of
$H^1(\R^n)$.

We next study the corresponding BMO space.
Since our atoms satisfy cancellation
against $\cP$, we define the BMO space
by using local approximation by elements
of $\cP$. More precisely, for a ball $B$,
let $\Pi_B$ be the orthogonal projection
in $L^2(B,\mu)$ onto $\cP|_B$, and define
\[
\|g\|_{\BMO(\mu)}
=\sup_{r_B\le1}\frac1{\mu(B)}
\int_B|g-\Pi_Bg|\,d\mu.
\]
We prove the John--Nirenberg inequality
and Hardy--BMO duality
\[
(H^1(\mu))^*\simeq\BMO(\mu)/\cP.
\]
See Theorems~\ref{thm:JN} and~\ref{thm:duality}.
We also prove the equivalence between
finite and infinite atomic norms and
establish interpolation theorems.

As a consequence of the Hardy--BMO duality,
we obtain
\[
\|R_jg\|_{\BMO(\mu)}
+\|R_j^*g\|_{\BMO(\mu)}
\lesssim\|g\|_{L^\infty(\mu)}.
\]
Together with the interpolation results,
this gives another proof of the $L^p(\mu)$
boundedness of the Riesz transforms for
$1<p<\infty$.

We finally consider the commutators of the
Riesz transforms. It is interesting to note
that $\BMO(\mu)$ is not the appropriate space
for the commutators. Indeed, when $n\ge2$,
we show that $x_2^2\in\BMO(\mu)$, but
$[x_2^2,R_2]$ is not bounded on $L^2(\mu)$.
See Proposition~\ref{prop:comm-counterexample}.

Instead, we consider the local BMO space
defined by
\[
\|b\|_{\BMO_{\mathrm{loc}}(\mu)}
=\sup_{r_B\le1}\frac1{\mu(B)}
\int_B|b-b_{B,\mu}|\,d\mu.
\]
We prove that, for every $1<p<\infty$
and $1\le j\le n$,
\[
\|[b,R_j]\|_{L^p(\mu)\to L^p(\mu)}
\simeq\|b\|_{\BMO_{\mathrm{loc}}(\mu)},
\]
with the same result for $R_j^*$.
See Theorem~\ref{thm:comm-characterization}.

In particular, we need only one Riesz
transform to characterize
$\BMO_{\mathrm{loc}}(\mu)$.
A similar result holds in the classical
setting~\cite{CRW}. However, for $n\ge2$,
the space $\BMO_{\mathrm{loc}}(\mu)$
is different from $\BMO(\mu)$.
Thus the BMO space associated with
the commutators is different from
the dual of $H^1(\mu)$.

The paper is organized as follows.
In Section~\ref{sec:spaces}, we study the
Hardy and BMO spaces, including atomic
decompositions, duality and interpolation.
Section~\ref{sec:Riesz} is devoted to
the Riesz transforms.
In Section~\ref{sec:commutators}, we study
the commutators of the Riesz transforms.
Finally, Section~\ref{sec:general-v} deals
with the general Laplacian with drift.

Throughout, $A\lesssim B$ means $A\le CB$,
where $C$ is independent of the variables
under consideration. We write $A\simeq B$
if both $A\lesssim B$ and $B\lesssim A$ hold.

\section{Atomic Hardy spaces and BMO spaces}\label{sec:spaces}

For every fixed $b>0$ and $0<r\le b$,
\begin{equation}\label{eq:local-volume}
 C_b^{-1}e^{2x_1}r^n
 \le \mu(B(x,r))\le C_b e^{2x_1}r^n.
\end{equation}
{Thus $\mu$ is locally doubling. More precisely, for every $b>0$, there exists $C_b>0$ such that, for $x\in\R^n$ and $\lambda>1$,}
\begin{equation}
	\label{eq-local doubling}
	\mu(B(x,\lambda r))\le C_b\lambda^n \mu(B(x,r)), \ \ \ 0<r<\lambda r<b.
\end{equation}

\begin{remark}
	\label{rem1}
	{The estimates \eqref{eq:local-volume} and \eqref{eq-local doubling} also hold for cubes in place of balls.}
\end{remark}

{The following space of polynomials will be used throughout the paper:}
\begin{equation*}
 \cP=\operatorname{span}_{\mathbb C}\{1,x_2,\ldots,x_n\}.
\end{equation*}
It is obvious that  $\cL p=0$ for all $p\in \cP$.

Let $E$ be a ball or cube. For
$f\in L^1(E,\mu)$, define $\Pi_Ef$ to be the unique element of
$\cP$ such that
\begin{equation}\label{eq:projection-definition}
 \int_E (f-\Pi_Ef)p\,d\mu=0,
 \qquad p\in\cP.
\end{equation}
{If $f\in L^2(E,\mu)$, then $\Pi_Ef$ is the unique minimizer of}
\[
 p\longmapsto \int_E|f-p|^2\,d\mu,
 \qquad p\in\cP.
\]

{To see that $\Pi_Ef$ exists and is unique, set} $p_0=1$ and $p_j=x_{j+1}$ for $1\le j\le n-1$. Assume that 
\[
 \Pi_Ef=\sum_{k=0}^{n-1}c_kp_k.
 \]
 From \eqref{eq:projection-definition}, we have
\[ 
 \sum_{k=0}^{n-1}
 \left(\int_Ep_jp_k\,d\mu\right)c_k
 =\int_Efp_j\,d\mu,
 \quad 0\le j\le n-1.
\]
We can rewrite the above equations as 
\begin{equation}\label{eq-matric eq}
A\textbf{c} = \textbf{b},
\end{equation}
where  $\textbf{c}=(c_0,c_1,\ldots, c_{n-1})^t$, {$\textbf{b}=\Big(\int_Efp_0\,d\mu,\ldots, \int_Efp_{n-1}\,d\mu\Big)^t$} and
\[
A =(a_{jk})_{j,k} =\Big(\int_Ep_jp_k\,d\mu\Big)_{j,k}.
\]
Since $A$ is positive definite, {i.e., $\xi^*A\xi>0$ for every $0\ne\xi\in\mathbb C^n$}, $A$ is non-singular. {It follows that} the equation \eqref{eq-matric eq} always has a unique solution $\textbf{c}$. {This ensures the existence} of the polynomial $\Pi_Ef$.

\begin{lemma}\label{lem-matrix}
	{Let $b>0$, and let $E$ be a ball (or cube) centered at $c$ with radius (or side length) $0<r\le b$. Define}
	\[
	p_0(x)=1,\qquad
	p_j(x)=\frac{x_{j+1}-c_{j+1}}r,
	\quad 1\le j\le n-1,
	\]
	and
	\[
	G_{jk}=\frac1{\mu(E)}\int_Ep_jp_k\,d\mu,
	\qquad 0\le j,k\le n-1.
	\]
	Then there exists $c_{b,n}>0$ such that
	\begin{equation*}
		\sum_{j,k=0}^{n-1}\overline{\xi_j}G_{jk}\xi_k
		\ge c_{b,n}\sum_{j=0}^{n-1}|\xi_j|^2
	\end{equation*}
	for all $\xi=(\xi_j)_j$; hence, $\|G^{-1}\|\le c_{b,n}^{-1}$.
\end{lemma}
\begin{proof}
{We give the proof for balls. The proof for cubes is similar.} 

Assume  $E=B(c,r)$.  It is easy to see that $|p_j(x)|\le 1$ for all $x\in E$ and $1\le j\le n-1$. Let
\[
G_{jk}=\frac1{\mu(E)}\int_Ep_jp_k\,d\mu,
\qquad 0\le j,k\le n-1.
\]
Using the change of variable 
$x=c+rz$, 
\[
G_{jk}=
\frac{\displaystyle\int_{E_0}\phi_j(z)\phi_k(z)e^{2rz_1}\,dz}
{\displaystyle\int_{E_0}e^{2rz_1}\,dz},
\quad \phi_0=1,\quad \phi_j(z)=z_{j+1},
\]
where $E_0=B(0,1)$ is a unit ball.

Note that $e^{2rz_1}\simeq_{b} 1$ and $\displaystyle\int_{E_0}e^{2rz_1}\,dz \simeq_{n,b} 1$. Hence,
\begin{equation*}
	\begin{aligned}
		\sum_{j,k=0}^{n-1}\overline{\xi_j}G_{jk}\xi_k
		&=\Big[ \int_{E_0}e^{2rz_1}\,dz\Big]^{-1 }\int_{E_0}
		{\left|\sum_{j=0}^{n-1}\xi_j\phi_j\right|^2e^{2rz_1}\,dz}\\
		&\gtrsim_{b,n} \int_{E_0}
		{\left|\sum_{j=0}^{n-1}\xi_j\phi_j\right|^2\,dz}.
	\end{aligned}
\end{equation*}

This, together with the fact that  $\phi_0,\ldots,\phi_{n-1}$ are linearly independent
on $E_0$, implies that  there exists $c_{b,n}>0$ such that
\begin{equation*}
	\sum_{j,k=0}^{n-1}\overline{\xi_j}G_{jk}\xi_k
	\ge c_{b,n}\sum_{j=0}^{n-1}|\xi_j|^2.
\end{equation*}
Thus $\|G^{-1}\|\le c_{b,n}^{-1}$ uniformly.

This completes our proof.	
\end{proof}
\begin{lemma}\label{lem:projection}
Let $b>0$ be fixed. If $E$ is a ball of radius at most $b$ or a cube of side
length at most $b$, then we have, for $1\le q\le\infty$,
\begin{equation}\label{eq:projection-Lq}
 \|\Pi_E f\|_{L^q(E,\mu)}\le C_{b,q}\|f\|_{L^q(E,\mu)}
\end{equation}
and
\begin{equation}\label{eq:projection-Linf}
 \|\Pi_E f\|_{L^\infty(E,\mu)}
 \le C_b\mu(E)^{-1}\|f\|_{L^1(E,\mu)}.
\end{equation}
Moreover,
\begin{equation}\label{eq:best-approx}
 \|f-\Pi_Ef\|_{L^q(E,\mu)}
 \simeq
 \inf_{p\in\cP}\|f-p\|_{L^q(E,\mu)},
\end{equation}
{with constants depending only on $b,n,q$.}
\end{lemma}

\begin{proof}
Consider  the basis
\[
p_0(x)=1,\qquad
p_j(x)=\frac{x_{j+1}-c_{j+1}}r,
\quad 1\le j\le n-1,
\]
for  $\cP$. Let $G =(G_{jk})$ be the matrix defined as in Lemma \ref{lem-matrix}. 

Write $\Pi_Ef=\sum_{k=0}^{n-1}a_kp_k$ and set
\[
 d_j=\frac1{\mu(E)}\int_Efp_j\,d\mu, \  \ \ 0\le j\le n-1.
\]
Then we have
\[
Ga = d,
\]
where $a=(a_k)_k$ and $d=(d_j)_j$.

 Since $|p_j|\le 1$ on $E$,
\[
 |d_j|\le \mu(E)^{-1}\|f\|_{L^1(E,\mu)},\qquad
 |a_k|\le C_{b,n}\mu(E)^{-1}\|f\|_{L^1(E,\mu)}.
\]
Thus
\begin{equation}\label{eq:projection-proof-Linfty}
 \|\Pi_Ef\|_{L^\infty(E,\mu)}
 \le C_{b,n}\mu(E)^{-1}\|f\|_{L^1(E,\mu)}.
\end{equation}
For $1\le q\le\infty$, using H\"older's inequality,
$$\|f\|_{L^1(E,\mu)}\le\mu(E)^{1-1/q}\|f\|_{L^q(E,\mu)}.
$$
Multiplying
\eqref{eq:projection-proof-Linfty} by $\mu(E)^{1/q}$, we have 
$${\|\Pi_Ef\|_{L^q(E,\mu)}\le C_{b,n}\|f\|_{L^q(E,\mu)}}.
$$
Finally, for every $p\in\cP$, $\Pi_Ep=p$. Therefore,
\[
 {\|f-\Pi_Ef\|_{L^q(E,\mu)}}
 ={\|(f-p)-\Pi_E(f-p)\|_{L^q(E,\mu)}}
 \le(1+C_{b,n}){\|f-p\|_{L^q(E,\mu)}}.
\]
Taking the infimum over $p\in\cP$ gives one inequality in
\eqref{eq:best-approx}; the reverse follows from $\Pi_Ef\in\cP$.

This completes our proof.
\end{proof}

{The next lemma is used to distribute cancellation across a finite family of overlapping balls or cubes.} We first introduce the following terms. Let $E_1,\ldots,E_N$ be subsets in $\mathbb R^n$. The overlap graph of
$\{E_i\}_{i=1}^N$ is the graph whose vertices
are $1,\ldots,N$, with an edge between $i$
and $j$ whenever $E_i\cap E_j\ne \emptyset$. A \textit{spanning tree} is a connected
subgraph containing all vertices and having
no cycles.

\begin{lemma}\label{lem:moment-transfer}
{Let $E_1,\ldots,E_N$ be balls (or cubes) with uniformly bounded and comparable radii (or side lengths).} Assume that the overlap graph is connected
and admits a spanning tree $T$ such that,
for every edge $(i,j)$ of $T$, the intersection
$E_i\cap E_j$ contains a ball $B_{ij}$ with
\[
r_{B_{ij}}\ge c\min\{ {\rm diam}E_i, {\rm diam}E_j\},
\]
where $c>0$ is independent of $i$ and $j$.

Set $E=\bigcup_{i=1}^N E_i$.
If $u\in L^q(E,\mu)$, $1<q\le\infty$, satisfies
\begin{equation*}
 \int up\,d\mu=0,\qquad p\in\cP,
\end{equation*}
then there are $u_i$ such that
\[
 u=\sum_{i=1}^N u_i,\qquad
 \supp u_i\subset E_i,\qquad
 \int u_i p\,d\mu=0  \ \ \ \text{for all $p\in\cP$}.
\]
In addition, for every $1<q\le\infty$,
\begin{equation}\label{eq:moment-transfer-bound}
 \mu(E_i)^{1-1/q}\|u_i\|_{L^q(E_i,\mu)}
 \le C\mu(E)^{1-1/q}\|u\|_{L^q(E,\mu)},\qquad 1\le i\le N.
\end{equation}
{The constant may depend on $N$, the bounds for the radii (or side lengths), and the overlap parameters, but not on $u$ or the location of the family.}
\end{lemma}

\begin{proof}

Choose a spanning tree $T$ as in the assumption.
Without loss of generality, we assume that $E_1$
is the root. For each vertex $i\ne1$, let $\pi(i)$
be its parent, that is, the next vertex on the path
from $i$ to the root. A vertex $j$ is called a child
of $i$ if $\pi(j)=i$, and a leaf is a vertex with
no children.

For each $k$, define
\[
\chi_k(x)=\frac{\mathbf1_{E_k}(x)}{\sum_{j=1}^N\mathbf1_{E_j}(x)},
\qquad x\in E,
\]
and set $\chi_k=0$ outside $E$. Put $f_k=\chi_k u$. Then
\[
u=\sum_{k=1}^N f_k,\qquad
\supp f_k\subset E_k,\qquad
\|f_k\|_{L^q(\mu)}\le\|u\|_{L^q(E,\mu)}.
\]

For each $2\le i\le N$, choose a ball
$D_i=B(c_i,r_i)\subset E_i\cap E_{\pi(i)}$
whose radius is comparable to the radii (or side lengths) of
$E_i$ and $E_{\pi(i)}$. Define
\[
p_{i,0}(x)=1,\qquad
p_{i,k}(x)=\frac{x_{k+1}-(c_i)_{k+1}}{r_i},
\quad 1\le k\le n-1,
\]
and
\[
A^{(i)}_{\ell m}=\int_{D_i}p_{i,\ell}p_{i,m}\,d\mu,
\qquad 0\le\ell,m\le n-1.
\]
By Lemma~\ref{lem-matrix}, $A^{(i)}$ is invertible. Set
\[
\theta_{i,k}=\mathbf1_{D_i}
\sum_{m=0}^{n-1}\bigl((A^{(i)})^{-1}\bigr)_{mk}p_{i,m},
\qquad 0\le k\le n-1.
\]
Then
\begin{equation}\label{eq-theta k ell}
\int\theta_{i,k}p_{i,\ell}\,d\mu=\delta_{k\ell},
\qquad 0\le k,\ell\le n-1.
\end{equation}

Let $i\ne1$ be any leaf of the current tree. Define
\[
h_i=\sum_{k=0}^{n-1}
\left(\int f_i p_{i,k}\,d\mu\right)\theta_{i,k},
\qquad u_i=f_i-h_i.
\]
Since $\supp h_i\subset D_i\subset E_i\cap E_{\pi(i)}$,
we have $\supp u_i\subset E_i$. Moreover,
by \eqref{eq-theta k ell},
\[
\int u_i p_{i,\ell}\,d\mu
=\int f_i p_{i,\ell}\,d\mu
-\sum_{k=0}^{n-1}
\left(\int f_i p_{i,k}\,d\mu\right)\delta_{k\ell}=0.
\]
Since $p_{i,0},\ldots,p_{i,n-1}$ form a basis of $\cP$, we have
\[
\int u_i p\,d\mu =0 \ \ \ \ \  \text{for all $p\in \cP$}.
\]
We now replace $f_{\pi(i)}$ by $f_{\pi(i)}+h_i$ and remove the leaf
$i$ from the tree.  {Repeat the same process with another leaf of the remaining tree. After $N-1$ steps, only the root remains.} Note that by Lemma~\ref{lem-matrix},
\[
\|\theta_{i,k}\|_{L^q(\mu)}
\le C\mu(D_i)^{1/q-1}.
\]
The overlap assumption and local doubling \eqref{eq-local doubling} imply
$\mu(D_i)\simeq\mu(E_i)$, while
$|p_{i,k}(x)|\le C$ for $x\in E_i$. Hence, by H\"older's
inequality,
\begin{align*}
\|h_i\|_{L^q(\mu)}
&\le C\mu(D_i)^{1/q-1}\|f_i\|_{L^1(\mu)}\\
&\le C\mu(D_i)^{1/q-1}\mu(E_i)^{1-1/q}
\|f_i\|_{L^q(\mu)}\\
&\le C\|f_i\|_{L^q(\mu)}.
\end{align*}
Thus $\|u_i\|_{L^q(\mu)}\le C\|f_i\|_{L^q(\mu)}$.
At each step, the only remaining function that changes is
$f_{\pi(i)}$, which is replaced by $f_{\pi(i)}+h_i$.
Since there are at most $N-1$ steps, the preceding estimates give
\[
\|u_i\|_{L^q(\mu)}\le C\|u\|_{L^q(E,\mu)},
\qquad 2\le i\le N,
\]
and the same bound for the function $f_1$ remaining at the root.

Finally, put $u_1=f_1$. We have
\[
u=\sum_{i=1}^N u_i,\qquad \supp u_1\subset E_1.
\]
Since $u$ and $u_2,\ldots,u_N$ have zero moments against $\cP$,
the same holds for $u_1$. As $\mu(E_i)\le\mu(E)$, we obtain
\[
\mu(E_i)^{1-1/q}\|u_i\|_{L^q(E_i,\mu)}
\le C\mu(E)^{1-1/q}\|u\|_{L^q(E,\mu)},
\qquad 1\le i\le N.
\]
This completes our proof.
\end{proof}

\subsection{Atomic Hardy spaces}\label{subsec:atoms}
\begin{definition}\label{def:H1}
Let $1<q\le\infty$ and $b>0$. A $(1,q;\cP)_b$-atom is a function $a$
supported in a ball $B=B(c,r)$ with $r\le b$ such that
\begin{equation*}
 \|a\|_{L^q(\mu)}\le \mu(B)^{1/q-1}
\end{equation*}
and
\begin{equation*}
 \int a(x)p(x)\,d\mu(x)=0,
 \qquad p\in\cP.
\end{equation*}
Equivalently,
\[
 \int a\,d\mu=0,
 \qquad
 \int x_j a(x)\,d\mu(x)=0,
 \quad 2\le j\le n.
\]
\end{definition}
\begin{definition}\label{def:H1 space} Let $1<q\le\infty$ and $b>0$. We define the Hardy space ${H^{1,q}_b(\mu)}$ as a set of all functions $f\in L^1(\mu)$ such that
\begin{equation}\label{eq:infinite-atomic-decomposition}
 f=\sum_{k=1}^\infty\lambda_ka_k\quad\text{in }L^1(\mu),
 \qquad \sum_{k=1}^\infty|\lambda_k|<\infty,
\end{equation}
where each $a_k$ is a $(1,q;\cP)_b$-atom. We also define 
\[
 \|f\|_{{H^{1,q}_b(\mu)}}
 =\inf\left\{\sum_{k=1}^\infty|\lambda_k|:
             f=\sum_{k=1}^\infty\lambda_ka_k\right\}.
\]
\end{definition}
Every $(1,q;\cP)_b$ atom has $L^1(\mu)$ norm at most $1$. Thus the series in
\eqref{eq:infinite-atomic-decomposition} converges absolutely in $L^1(\mu)$,
and
\begin{equation*}
 \|f\|_{L^1(\mu)}\le \|f\|_{{H^{1,q}_b(\mu)}}.
\end{equation*}

\begin{lemma}\label{lem-completeness}
	For $1<q\le \infty$ and $b>0$, the Hardy space ${H^{1,q}_b(\mu)}$ is a Banach space.
\end{lemma}
\begin{proof}
	The proof is quite {standard}, but for the sake of completeness, we will provide it here. 
	
	It suffices to prove ${H^{1,q}_b(\mu)}$ is complete.  {Let $\{f_j\}$ be a Cauchy sequence} in ${H^{1,q}_b(\mu)}$. We may assume that
	\[
	\|f_{j+1}-f_j\|_{{H^{1,q}_b(\mu)}}\le 2^{-j},
	\qquad j\ge 1;
	\]
	otherwise, we consider a subsequence.

	Set $g_1=f_1$ and $g_j=f_j-f_{j-1}$ for $j\ge 2$. Then
	\[
	\sum_{j=1}^\infty \|g_j\|_{{H^{1,q}_b(\mu)}}<\infty.
	\]
	For each $j$, choose an atomic decomposition
	\[
	g_j=\sum_{k=1}^\infty \lambda_{j,k}a_{j,k}
	\]
	such that
	\[
	\sum_{k=1}^\infty|\lambda_{j,k}|
	\le \|g_j\|_{{H^{1,q}_b(\mu)}}+2^{-j}.
	\]
	Hence
	\[
	\sum_{j,k}|\lambda_{j,k}|<\infty.
	\]
	Since $\|a_{j,k}\|_{L^1(\mu)}\le 1$ for each $j,k$,
	the series
	\[
	f=\sum_{j,k}\lambda_{j,k}a_{j,k}
	\]
	converges in $L^1(\mu)$ and defines an element of ${H^{1,q}_b(\mu)}$.
	Moreover,
	\[
	\|f-f_m\|_{{H^{1,q}_b(\mu)}}
	\le \sum_{j=m+1}^\infty
	\left(\|g_j\|_{{H^{1,q}_b(\mu)}}+2^{-j}\right)
	\longrightarrow 0.
	\]
	Thus ${H^{1,q}_b(\mu)}$ is complete.
	
	This completes our {proof}.

\end{proof}

\begin{proposition}\label{prop:H1-independence}
For $1<q\le\infty$ and $b>0$, we have
$$
{H^{1,q}_b(\mu)} = {H^{1,\infty}_1(\mu)} 
$$
 with equivalent norms. Therefore, we will define ${H^1(\mu)}$ as any space ${H^{1,q}_b(\mu)}$ with $1<q\le \infty$ and $b>0$. 
\end{proposition}

\begin{proof}
 Let $0<b_1<b_2$. If $a$ is  an $(1,q;\cP)_{b_1}$-atom, then  $a$ is  also an {$(1,q;\cP)_{b_2}$-atom}. Hence, ${H^{1,q}_{b_1}(\mu)}\hookrightarrow {H^{1,q}_{b_2}(\mu)}$.

 Conversely, if $a$ is  an $(1,q;\cP)_{b_2}$-atom associated with $B$ with $b_1<r_B\le b_2$, then {we can cover $B$ by a family of balls} $\{B_i: i \in \mathcal I\}$ satisfying the following properties:
 \begin{enumerate}
 	\item[(i)] $r_{B_i}=b_1$ and $B_i\subset 2B$ for all $i\in \mathcal I$;
 	\item[(ii)]  {the overlap graph admits a spanning tree $T$ such that $B_i\cap B_j$ contains a ball of radius $b_1/10$ for every edge $(i,j)$ of $T$;}
 	\item[(iii)] $\sharp \mathcal I \lesssim (b_2/b_1)^n$;
 	 \item[(iv)] $\sum_{i\in \mathcal I}\mathbf{1}_{B_i}\lesssim 1$.
 \end{enumerate}
By \eqref{eq:moment-transfer-bound}, we can write 
\[
a = \sum_{i\in \mathcal I} u_i,
\]
where $u_i$ satisfies
\[
\supp u_i\subset B_i,\qquad
\int u_i p\,d\mu=0  \ \ \ \text{for all $p\in\cP$},
\]
and
\[
 \mu(B_i)^{1-1/q}\|u_i\|_{L^q(\mu)}
\lesssim \mu(2B)^{1-1/q}\|{a}\|_{L^q(\mu)}\lesssim \mu(2B)^{1-1/q}\mu(B)^{1/q-1}\lesssim 1.
\]
This implies
\[
\|u_i\|_{L^q(\mu)}\lesssim  \mu(B_i)^{1/q-1}.
\]
That is, $u_i$ is a multiple of an $(1,q;\cP)_{b_1}$ atom. Consequently, $a\in {H^{1,q}_{b_1}(\mu)}$ with
\[
\|a\|_{{H^{1,q}_{b_1}(\mu)}}\lesssim_{b_1,b_2,n} 1.
\]
This implies, for $1<q\le \infty$ 
\[
 {H^{1,q}_{b_1}(\mu)}= {H^{1,q}_{b_2}(\mu)}  \ \  \text{for all $0<b_1<b_2$}.
\]

To complete the proof, we need to show that  $ {H^{1,q}_{1}(\mu)}= {H^{1,\infty}_{1}(\mu)}$ for $1<q<\infty$.  Fix $1<q<\infty$. Since \[
{H^{1,q}_{b_1}(\mu)}= {H^{1,q}_{b_2}(\mu)}  \ \  \text{for all $0<b_1<b_2$},
\]
it suffices to prove that if $a$ is a  {cube $(1,q;\cP)$-atom} supported in a cube $Q$ with {$\ell(Q)\le1$} and $\|a\|_{L^q(\mu)}\le \mu(Q)^{1/q-1}$, then $a\in {H^{1,\infty}_{1}(\mu)}$.  Let $K>0$ be positive large number which will be fixed later. Denote by $\{Q_i\}$ the family of  dyadic subcubes $Q_i\subset Q$ satisfying
\[
 \frac1{\mu(Q_i)}\int_{Q_i}|a|^q\,d\mu
   >\left(\frac K{\mu(Q)}\right)^q.
\]
The local doubling property \eqref{eq-local doubling} and the maximality of $Q_i$ imply
\begin{equation}\label{eq-sum Qi}
 \frac1{\mu(Q_i)}\int_{Q_i}|a|^q\,d\mu
 \le C\left(\frac K{\mu(Q)}\right)^q,
 \qquad
 \sum_i\mu(Q_i)\le K^{-q}\mu(Q).
\end{equation}
Set
\[
 b_i=(a-\Pi_{Q_i}a)\mathbf1_{Q_i},
 \qquad
 g=a-\sum_i b_i.
\]
From \eqref{eq:projection-definition},  $b_i$ annihilates $\cP$, and so does $g$.  {That is,}
\[
\int b_i p \,d\mu = \int g p \,d\mu = 0 \ \ \ \text{for all $p\in \cP$}. 
\]
We now estimate $\|g\|_{\infty}$. From the definition of $g$,
\[
g =a \mathbf{1}_{Q\setminus \cup_{i}Q_i} +\sum_{i} (\Pi_{Q_i}a) \mathbf1_{Q_i}
\]

 By \eqref{eq-local doubling} and dyadic Lebesgue differentiation theorem,
 \[
 \Big\|a \mathbf{1}_{Q\setminus \cup_{i}Q_i}\Big\|_\infty \le \frac{CK}{\mu(Q)}.
 \]
In addition, by   \eqref{eq:projection-Linf}, H\"older's inequality and \eqref{eq-sum Qi},
\[
\begin{aligned}
	\Big\|\sum_{i} (\Pi_{Q_i}a) \mathbf1_{Q_i}\Big\|_\infty &{\le \sup_{i} \|\Pi_{Q_i}a\|_\infty}\\
	&{\le A_0 \sup_i \mu(Q_i)^{-1}\|a\|_{L^1(Q_i,\mu)}}\\
	&{\le  A_0\sup_i \Big(\frac{1}{\mu(Q_i)}\int_{Q_i} |a|^q\,d\mu\Big)^{1/q}}\\
	&\le \frac{A_0K}{\mu(Q)},
\end{aligned}
\]
where $A_0$ depends on $n$ only.

Consequently,
\[
 \|g\|_\infty\le \frac{CK}{\mu(Q)}.
\]
Thus $g$ is $CK$ times an  $(1,\infty; \cP)_1$ atom associated with the cube $Q$. Moreover, by H\"older's {inequality}, \eqref{eq:projection-Lq} and \eqref{eq-sum Qi},
\begin{align*}
 \sum_i \mu(Q_i)^{1-1/q}\|b_i\|_{L^q(\mu)}
 & {\le } C_0 \left(\sum_i\mu(Q_i)\right)^{1-1/q}
       \left(\sum_i\|a\|_{L^q(Q_i,\mu)}^q\right)^{1/q}\\
       &\le  C_0\left(K^{-q}\mu(Q)\right)^{1-1/q}
        \|a\|_{L^q(Q ,\mu)}\\
 & {\le} C_0K^{-(q-1)}.
\end{align*}
 Choose
$K$ so large that {$C_0K^{-(q-1)}\le 1/2$}. Then we have
\[
 {\sum_i \mu(Q_i)^{1-1/q}\|b_i\|_{L^q(\mu)}\le \frac{1}{2}.}
\]
Hence, by normalization, we can write 
\[
{a=g+\sum_i\lambda_i a_i=:g_1+f_1,}
\]
{where $g_1=g$ is a multiple of an $(1,\infty;\cP)_1$-atom associated with $Q$, each $a_i$}
is an $(1,q;\mathcal P)_1$-atom supported in the dyadic subcube
$Q_i$, and 
\[
{\|g_1\|_{H^{1,\infty}_1(\mu)}\le A_0K,\qquad \sum_i|\lambda_i|\le\frac12.}
\]

{Repeating the construction for each $a_i$, we obtain} 
\[
  a_{ i}=\beta_i g_{i}
   +\sum_j\lambda_{i,j}a_{i,j},
\]
where $g_i$ is an $(1,\infty,\cP)_1$ atom for each $i$, $a_{i,j}$
	is an $(1,q;\mathcal P)_1$-atom for each $i,j$ and
	\[
	|\beta_i|\le  A_0K ,\ \ \ { \sum_j|\lambda_{i,j}|\le\frac12.}
	\]

Hence, 
\[
f_1 =   \sum_{i}\lambda_i \beta_i g_i + \sum_{i,j}\lambda_i \lambda_{i,j} a_{i,j}=:g_2 + f_2.
\]
where $ g_i$ is an $(1,\infty,\cP)_1$ atoms for each $i$, $a_{i,j}$
is an $(1,q;\mathcal P)_1$-atom for each $i,j$; in addition,
\[
{\|g_2\|_{{H^{1,\infty}_1(\mu)}}}\le  \sum_{i}|\lambda_i \beta_i|\le   \frac{1}{2}A_0K, \ \ \ \ \text{and}\ \   {\|f_2\|_{H^{1,q}_1(\mu)}\le \sum_{i,j}|\lambda_i \lambda_{i,j}|\le \frac{1}{2^2}.}
\]

By iteration, we have
\[
\begin{aligned}
	a &= f_1 + g_1 \\
	f_1 &= f_2 + g_2 \\
	f_2 &= f_3 + {g_3} \\
	\vdots&\ \hskip 1cm\vdots 
\end{aligned}
\]
where 
\begin{equation}\label{eq- gk norm}
{\|g_k\|_{{H^{1,\infty}_1(\mu)}}\le\frac{A_0K}{2^{k-1}}}, \ \ \ \ \text{and}\ \   {\|f_k\|_{{H^{1,q}_1(\mu)}}\le\frac{1}{2^k}}.
\end{equation}
Hence, for each $k$, 
\[
a= g_1 +\ldots + g_k + f_k:= \textbf{g}_k + f_k.
\]
{By \eqref{eq- gk norm} and Lemma~\ref{lem-completeness}}, $\lim \textbf{g}_k= {\textbf{g}\in {H^{1,\infty}_1(\mu)}}$. In addition, by Fatou's lemma,
\[
\|a-\textbf{g}\|_{L^1(\mu)} \le \lim \|a-\textbf{g}_k\|_{L^1(\mu)} =\lim \|f_k\|_{L^1(\mu)}\le \lim {\|f_k\|_{H^{1,q}_1(\mu)}} \le \lim \frac{1}{2^k} =0.
\] 
This implies that ${a=\textbf{g} \in {H^{1,\infty}_1(\mu)}}$ and this completes our proof.
\end{proof}

\begin{definition}\label{def:BMO}
{Let $1\le q<\infty$ and $b>0$. The space ${\BMO_b^q(\mu)}$ is defined} as a set of all $f\in L^1_{\rm loc}(\mu)$ such that 
\begin{equation*}
 \|f\|_{{\BMO_b^q(\mu)}}
{:=\sup_{B:\,r_B\le b}}
 \left(\frac1{\mu(B)}\int_B|f-\Pi_Bf|^q\,d\mu\right)^{1/q}<\infty.
\end{equation*}
\end{definition}

\begin{proposition}\label{prop:BMO-scale}
Let $b_1,b_2>0$. The seminorms $\|\cdot\|_{{\BMO_{b_1}^q(\mu)}}$ and
$\|\cdot\|_{{\BMO_{b_2}^q(\mu)}}$ are equivalent. The same space and an
equivalent seminorm are obtained if balls are replaced by cubes in Definition \ref{def:BMO}.
\end{proposition}

\begin{proof}
Let $b_1<b_2$. Obviously, $\|f\|_{{\BMO_{b_1}^q(\mu)}}\le \|f\|_{{\BMO_{b_2}^q(\mu)}}$. We now prove 
\[
\|f\|_{{\BMO_{b_2}^q(\mu)}}\lesssim \|f\|_{{\BMO_{b_1}^q(\mu)}}.
\]
Let $B$ be a ball with $b_1<r_B\le b_2$. Then {we can cover $B$ by a family of balls} $\{B_i: i \in \mathcal I\}$ satisfying the following properties:
\begin{enumerate}
	\item[(i)] $r_{B_i}=b_1$ and $B_i\subset 2B$ for all $i\in \mathcal I$;
	\item[(ii)]  {the overlap graph admits a spanning tree $T$ such that $B_i\cap B_j$ contains a ball of radius $b_1/10$ for every edge $(i,j)$ of $T$;}
	\item[(iii)] $\sharp \mathcal I \lesssim (b_2/b_1)^n$;
	\item[(iv)] $\sum_{i\in \mathcal I}\mathbf{1}_{B_i}\lesssim 1$.
\end{enumerate}

Set $p_i=\Pi_{B_i}f$. For two balls $B_i,B_j$ with $(i,j)\in T$, let $D_{ij}\subset B_i\cap B_j$ be a ball with radius $b_1/10$. We have
\[
 p_i-p_j=\Pi_{D_{ij}}(p_i-f)-\Pi_{D_{ij}}(p_j-f).
\]
This, together with  Lemma~\ref{lem:projection} and \eqref{eq-local doubling}, implies
\[
\begin{aligned}
	\|p_i-p_j\|_{L^\infty(D_{ij},\mu)}&\le \|\Pi_{D_{ij}}(p_i-f)\|_{L^\infty(D_{ij},\mu)}+\|\Pi_{D_{ij}}(p_j-f)\|_{L^\infty(D_{ij},\mu)}\\
	&\lesssim \mu(D_{ij})^{-1}\|p_i-f \|_{L^1(D_{ij},\mu)}+\mu(D_{ij})^{-1}\| p_j-f \|_{L^1(D_{ij},\mu)}\\
	&\lesssim \frac1{\mu(B_i)}\int_{B_i}|f-\Pi_{B_i}f|\,d\mu +\frac1{\mu(B_j)}\int_{B_j}|f-\Pi_{B_j}f|\,d\mu\\
	&\lesssim \|f\|_{{\BMO_{b_1}^q(\mu)}}.
\end{aligned}
\]
{Since the family of balls is connected}, the above inequality implies
\[
	\|p_i-p_j\|_{L^\infty(D_{ij},\mu)}\lesssim_{b_1,b_2} {\|f\|_{{\BMO_{b_1}^q(\mu)}}}\ \ \ \  \text{for all $i,j \in \mathcal I$}.
\]
This, together with the local doubling \eqref{eq-local doubling}, further yields
\[
\|p_i-p_j\|_{L^\infty(B,\mu)}\lesssim_{b_1,b_2} {\|f\|_{{\BMO_{b_1}^q(\mu)}}} \ \ \ \  \text{for all $i,j \in \mathcal I$}.
\]

It follows that 
$$\|p_i-\bar p\|_{L^\infty(B,\mu)}\lesssim_{b_1,b_2}  \|f\|_{{\BMO_{b_1}^q(\mu)}}
$$
where {$\bar p=|\mathcal I|^{-1}\sum_{i\in\mathcal I}p_i\in\cP$.}

Consequently,
\[
 \begin{aligned}
 \int_B|f-\bar p|^q\,d\mu&\lesssim_{b_1,b_2}\sum_i\bigl(\|f-p_i\|_{L^q(B_i,\mu)}^q
 +\mu(B_i)\|p_i-\bar p\|_{L^\infty(B_i,\mu)}^q \bigr)\\
 &\lesssim_{b_1,b_2}\sum_i\bigl(\|f-p_i\|_{L^q(B_i,\mu)}^q
                  +\mu(B_i)\|f\|^q_{{\BMO_{b_1}^q(\mu)}}\bigr)\\
 &\le C\mu(B)\|f\|^q_{{\BMO_{b_1}^q(\mu)}}.
 \end{aligned}
\]
This, in combination with  \eqref{eq:best-approx}, implies
\[
\sup_{B:b_1<\,r_B\le b_2}
\left(\frac1{\mu(B)}\int_B|f-\Pi_Bf|^q\,d\mu\right)^{1/q}\lesssim \|f\|_{{\BMO_{b_1}^q(\mu)}}.
\]
This completes our proof.

\end{proof}

\subsection{The John--Nirenberg inequality}\label{sec:JN}

\begin{theorem}\label{thm:JN}
There are constants $C,c>0$ such that, for every $f\in{\BMO_1^1(\mu)}$, every
cube $Q$ with $\ell(Q)\le 1$ and every $\lambda>0$,
\begin{equation*}
 \mu\{x\in Q:|f(x)-\Pi_Qf(x)|>\lambda\}
 \le C\exp\left(-c\lambda/{\|f\|_{{\BMO_1^1(\mu)}}}\right)\mu(Q).
\end{equation*}
 Consequently,
for every $1\le q<\infty$,
\begin{equation*}
 \|f\|_{{\BMO_1^q(\mu)}}\simeq_q \|f\|_{{\BMO_1^1(\mu)}}.
\end{equation*}
If $\|f\|_{{\BMO_1^1(\mu)}}=0$, then $f\in\cP$ almost everywhere.
\end{theorem}

\begin{proof}
{By Proposition~\ref{prop:BMO-scale}, we may work with cubes of side length at most $1$.}
Let $Q$ be  a cube with $\ell(Q)\le 1$. Define
\[
 \beta=\|f\|_{{\BMO_1^1(\mu)}}.
\]
{If $\beta=0$, then $f=\Pi_Bf$ almost everywhere on every ball $B$ with $r_B\le1$. Since these polynomials agree on overlapping balls, $f\in\cP$ and the result is immediate. We may therefore assume that $\beta>0$.}
{By Lemma~\ref{lem:projection} and Remark~\ref{rem1}, we have}
\begin{equation*}
 \frac1{\mu(S)}\int_S|f-\Pi_Sf|\,d\mu\le A\beta
\end{equation*}
for every dyadic subcube $S\subset Q$, where $A$ is independent of $S$.

Choose $T>2$. Let $\{Q_i\}_i$ be the set of maximal dyadic subcubes in $Q$ satisfying
\begin{equation}\label{eq-stop cond}
 \frac1{\mu(Q_i)}\int_{Q_i}|f-\Pi_Qf|\,d\mu>TA\beta.
\end{equation}
Then the family $\{Q_i\}_i$ is pairwise disjoint and
\begin{equation}\label{eq:first-bad-measure}
\begin{aligned}
	\sum_i\mu(Q_i)&\le \frac{1}{TA\beta} \sum_{i}\int_{Q_i}|f-\Pi_Qf|\,d\mu\\
	&\le \frac{1}{TA\beta} \int_{Q}|f-\Pi_Qf|\,d\mu\\
	&\le T^{-1}\mu(Q).
\end{aligned} 
\end{equation}
For each $Q_i$, let $R_i$ be its parent. Since  $R_i$ does not satisfy the stopping condition \eqref{eq-stop cond}, we have
\begin{equation}\label{eq:stopped-average-upper}
\begin{aligned}
	 \frac1{\mu(Q_i)}\int_{Q_i}|f-\Pi_Qf|\,d\mu&\le \frac{\mu(R_i)}{\mu(Q_i)}\frac1{\mu(R_i)}\int_{R_i}|f-\Pi_Qf|\,d\mu\\
	&\le C\frac1{\mu(R_i)}\int_{R_i}|f-\Pi_Qf|\,d\mu\\
	&\le C T A\beta,
\end{aligned}	
\end{equation}
where in the second inequality we use the local doubling \eqref{eq-local doubling} for cubes.

Hence, this, {in combination} with  $\Pi_{Q_i}\Pi_Qf=\Pi_Qf$ and
\eqref{eq:projection-Linf}, implies
\begin{equation}\label{eq:poly-jump}
	\begin{aligned}
		\|\Pi_{Q_i}f-\Pi_Qf\|_{L^\infty(Q_i,\mu)}&=\| \Pi_{Q_i}(f-\Pi_Qf)\|_{L^\infty(Q_i,\mu)}\\
		&\le C \mu(Q_i)^{-1}\| f-\Pi_Qf \|_{L^1(Q_i,\mu)}\\
		&\le C T A\beta.
	\end{aligned}
\end{equation}
In addition, the dyadic Lebesgue differentiation theorem implies
\[
|f-\Pi_Qf|\le TA\beta, \ \ \ \ \ \text{a.e. on $Q\setminus  \bigcup_iQ_i$}.
\]

Repeat the same construction for every $Q_i$ with
$f-\Pi_{Q_i}f$, we can find a family of subcubes $\{Q_{i,j}\}_{i,j}$ such that for each $i$, 
\begin{equation}\label{eq-inequality for Qj}
\sum_{j}\mu(Q_{i,j})\le T^{-1}\mu(Q_i)
\end{equation}
and
\begin{equation}\label{eq-inequality for i}
|f-\Pi_{Q_i}f|\le TA\beta, \ \ \ \ \ \text{a.e. on $Q_i\setminus  \bigcup_jQ_{i,j}$}.
\end{equation}
It follows from \eqref{eq-inequality for Qj} and \eqref{eq:first-bad-measure} that
\[
\sum_{i,j}\mu(Q_{i,j})\le T^{-1}\sum_i \mu(Q_i)\le T^{-2}\mu(Q).
\]
If $x\in Q\setminus \bigcup_iQ_{i}$, then
\[
|f-\Pi_Qf|\le TA\beta.
\]
If $x\in \big(\bigcup_iQ_{i}\big)\setminus \big( \bigcup_{i,j}Q_{i,j}\big)=\cup_i \Big(Q_i\setminus \bigcup_jQ_{i,j}\Big)$, then  $x\in  Q_i\setminus \bigcup_jQ_{i,j}$ for some $i$. For such an $i$,
using \eqref{eq-inequality for i} and \eqref{eq:poly-jump},
\[
|f-\Pi_{Q}f|\le |f-\Pi_{Q_i}f| + |\Pi_{Q}f-\Pi_{Q_i}f|\le TA\beta +CTA\beta\le 2CTA\beta  
\]
for $x\in  Q_i\setminus \bigcup_jQ_{i,j}$.

Consequently,
\[
|f-\Pi_{Q}f|\le 2CTA\beta  \ \ \ \ \ \text{a.e. on $Q\setminus  \bigcup_{i,j}Q_{i,j}$}.
\]
Repeating the same construction, after $k$ generations we will obtain a family of cubes which is denoted, without any confusion, by $\{Q^k_j\}$ such that 
\begin{equation*}
	\sum_{i}\mu(Q^k_{i})\le T^{-k}\mu(Q)
\end{equation*}
and
\begin{equation*}
	|f-\Pi_{Q}f|\le CkTA\beta, \ \ \ \ \ \text{a.e. on $Q\setminus  \bigcup_iQ^k_{i}$}.
\end{equation*}

Therefore
\[
 \mu\{x\in Q:|f-\Pi_Qf|>CkTA\beta\}
 \le T^{-k}\mu(Q).
\]
For $\lambda\ge L\beta$, take
$k=\lfloor\lambda/(L\beta)\rfloor$. For $0<\lambda<L\beta$,
use the trivial bound by $\mu(Q)$. This gives, for every $\lambda>0$,
\[
\mu\{x\in Q:|f(x)-\Pi_Qf(x)|>\lambda\}
\le T\exp\left(-\frac{(\log T)\lambda}{L\beta}\right)\mu(Q).
\]

This completes our proof.

\end{proof}

{We write $\BMO(\mu)=\BMO_1^1(\mu)$. By Theorem~\ref{thm:JN}, the null space of this seminorm is $\cP$.}

\subsection{Hardy--BMO duality}\label{sec:duality}

\begin{lemma}\label{lem:bounded-approximation}
Let $g\in{\BMO(\mu)}$. There are $p\in\cP$ and functions
$g_m\in L^\infty(\mu)$ such that
\begin{equation*}
 \sup_m\|g_m\|_{{\BMO(\mu)}}\le C\|g\|_{{\BMO(\mu)}},
 \qquad g_m\longrightarrow g-p
 \quad\text{in }L^s_{\mathrm{loc}}(\mu),\quad 1\le s<\infty.
\end{equation*}
\end{lemma}

\begin{proof}
Set $\beta=\|g\|_{{\BMO(\mu)}}$. If $\beta=0$, then $g\in\cP$ by
Theorem~\ref{thm:JN}; take $p=g$ and $g_m=0$. Assume $\beta>0$.
All convolutions used below are with respect to Lebesgue measure. Choose a
nonnegative radial function $\eta\in C_c^\infty(B(0,1/8))$ with $\displaystyle \int\eta\,dx=1$,
and put $G=\eta*g$ and $h=g-G$.


Since $\eta$ is radial, $\displaystyle \int z_j\eta(z)\,dz=0$ for every $j$.
Thus, for $q\in\cP$,
\[
 (\eta*q)(x)=\int\eta(z)q(x-z)\,dz=q(x).
\]
To prove the local estimate for $h$, set $p_c=\Pi_{B(c,2)}g$.
Since $p_c\in \cP$,  using the above identity, we have
\[
 h=(g-p_c)-\eta*(g-p_c).
\]
Using $d\mu(x)=e^{2z_1}d\mu(x-z)$,
\[
\begin{aligned}
 \int_{B(c,1)}|\eta*(g-p_c)|^2\,d\mu
 &\le\int_{|z|<1/8}\eta(z)e^{2z_1}
       \int_{B(c,1)-z}|g-p_c|^2\,d\mu\,dz\\
 &\le C\int_{B(c,2)}|g-p_c|^2\,d\mu\\
 &\le C{\beta^2\mu(B(c,2))},
\end{aligned}
\]
where in the last inequality we used 
Theorem~\ref{thm:JN}.

This implies, for  $c\in \mathbb R^n$,
\begin{equation}\label{eq:approx-local-remainder}
 \begin{aligned}
 	   \left(\frac1{\mu(B(c,1))}\int_{B(c,1)}|h|^2\,d\mu\right)^{1/2}&\le  \left(\frac1{\mu(B(c,1))}\int_{B(c,1)}|g-p_c|^2\,d\mu\right)^{1/2}\\
 	   & \ \ \ \ \ +\left(\frac1{\mu(B(c,1))}\int_{B(c,1)}|\eta\ast (g-p_c)|^2\,d\mu\right)^{1/2}\\
 	   &\le C\beta.
 \end{aligned}
\end{equation}

Next, we note that for any $p\in \cP$,
$$
(\partial_1\eta)*p=\partial_1p=0, \ \ \ \ (\partial_j\partial_k\eta)*p=\partial_j\partial_kp=0, \ \ \ \ j,k\ge 2.
$$
For $x\in\R^n$, set $p_x=\Pi_{B(x,1)}g$. If $K\in \{\partial_1\eta, \partial_j\partial_k\eta: j,k\ge 2\}$, then
\[
\begin{aligned}
	|(K*g)(x)|
	&=\left|\int K(z)(g-p_x)(x-z)\,dz\right|\\
	&\le C\int_{B(x,1)}|g-p_x|\,dy\\
	&\le C\int_{B(x,1)}|g-p_x|\,e^{-2y_1}e^{2y_1}dy.
\end{aligned}
\]
By \eqref{eq:local-volume}, we have 
$$e^{-2y_1} ={e^{-2x_1} e^{2(x_1-y_1)}}\simeq  {e^{-2x_1} \simeq \mu(B(x,1))^{-1}}, \qquad \qquad y\in B(x,1).
$$ 
This further implies
\[
\begin{aligned}
	|(K*g)(x)|	&\le \frac{C}{\mu(B(x,1))}\int_{B(x,1)}|g-p_x|\,d\mu(y)\\
	&\le C\beta.
\end{aligned}
\]
It follows that
\begin{equation}\label{eq:approx-smooth-derivatives}
 \|\partial_1G\|_\infty+
 \sum_{j,k=2}^n\|\partial_j\partial_kG\|_\infty\le C\beta.
\end{equation}

Write $x=(t,y)\in\R\times\R^{n-1}$ and set
\[
 p(y)=G(0,0)+\nabla_yG(0,0)\cdot y,\qquad F(t,y)=G(t,y)-p(y).
\]
Using \eqref{eq:approx-smooth-derivatives} and the Taylor expansion, we have 
$$|F(t,y)-F(0,y)|\le C\beta|t|, \ \ \ \ \ \ \ \ \ \ \ \ \ |F(0,y)|\le C\beta|y|^2.
$$
 Hence
$|F(t,y)|\le C\beta(|t|+|y|^2)$.
To estimate $\nabla_yF$, set $H(y)=F(t,y)-F(0,y)$.
Then $\|H\|_\infty\le C\beta|t|$. In addition, {since $D_y^2p=0$}, we have, by \eqref{eq:approx-smooth-derivatives},
\[
\|D_y^2H\|_\infty\le \|D_y^2G(t,y)\|_\infty+\|D_y^2G(0,y)\|_\infty\le C\beta.
\]
{Taylor's formula gives} 
\[
 |\partial_jH(y)|\le C\left(\frac{\|H\|_\infty}{s}
                         +s\|\partial_j^2H\|_\infty\right),
 \qquad s>0.
\]
Taking $s=\sqrt{|t|}$, we have
$|\nabla_yH|\le C\beta\sqrt{|t|}$.
Also, since
\[
F(0,y)=G(0,y) -G(0,0) -{\nabla_yG(0,0)\cdot y}
\]
{it follows from the mean value theorem and \eqref{eq:approx-smooth-derivatives} that}
\[
|\nabla _yF(0,y)|=|\nabla_yG(0,y) - \nabla_yG(0,0)|\le C|y|\beta.
\]
Thus
\begin{equation}\label{eq:approx-growth}
 |F(t,y)|\le C\beta(|t|+|y|^2),\qquad
 |\nabla_yF(t,y)|\le C\beta(|y|+\sqrt{|t|}).
\end{equation}

Choose compactly supported smooth cutoffs $\chi_1$ on $\R$ and
$\chi_\perp$ on $\R^{n-1}$, equal to one near the origin.
For $R\ge1$, put
\[
 F_R(t,y)=\chi_1(t/R^2)\chi_\perp(y/R)F(t,y).
\]
On the support of the cutoff functions, $|t|\le CR^2$, $|y|\le CR$.
Therefore \eqref{eq:approx-growth} implies
$$
|F_R|\le C\beta R^2, \ \ \ \ \ |\nabla_yF_R|\le C\beta R
$$
The product rule, together with
$|\partial_t\chi_1(t/R^2)|\le CR^{-2}$,
$|\nabla_y\chi_\perp(y/R)|\le CR^{-1}$,
$|D_y^2\chi_\perp(y/R)|\le CR^{-2}$ and \eqref{eq:approx-smooth-derivatives}, yields
\begin{equation}\label{eq:approx-cutoff-derivatives}
 \|\partial_1F_R\|_\infty+
 \sum_{j,k=2}^n\|\partial_j\partial_kF_R\|_\infty
 \le C\beta,
\end{equation}
uniformly in $R$.

For a ball $B(c,r)$, $r\le1$,  {set $c_\perp=(c_2,\ldots,c_n)\in\R^{n-1}$}.
{By Taylor's expansion and \eqref{eq:approx-cutoff-derivatives},}
\[
 \sup_{x=(t,y)\in B(c,r)}
 \big|F_R(x)-F_R(c)-\nabla_yF_R(c)\cdot(y-c_\perp)\big|
 \le C\beta(r+r^2)\le C\beta.
\]
{This, together with the fact that $F_R(c)+\nabla_yF_R(c)\cdot(y-c_\perp)\in\cP$} and \eqref{eq:best-approx}, implies $\|F_R\|_{{\BMO(\mu)}}\le C\beta$.
Similarly, we also obtain
$\|G\|_{{\BMO(\mu)}}\le C\beta$. Hence
$\|h\|_{{\BMO(\mu)}}\le C\beta$ since $h=g-G$.

For $0<\varepsilon\le1$, set
$\eta_\varepsilon(x)=\varepsilon^{-n}\eta(x/\varepsilon)$. We have, for $B=B(c,r)$ with $r<1$,
\[
\begin{aligned}
	\frac{1}{\mu(B)}{\int_B} |\eta_\varepsilon*h(x) &- \Pi_B(\eta_\varepsilon*h){(x)}|\,d\mu(x)\\
	&=\frac{1}{\mu(B)}\int_B\Big|\int_{|z|<\varepsilon/8}{\eta_\varepsilon}(z)\Big[h(x-z) -(\Pi_Bh(\cdot-z))(x)\Big]\, dz\Big| d\mu(x).
\end{aligned}
\]
Using the fact $\Pi_B[u(\cdot-z)](x)=(\Pi_{B-z}u)(x-z)$, {we further obtain}
\[
\begin{aligned}
	\frac{1}{\mu(B)}\int_B |\eta_\varepsilon*h(x) &- \Pi_B(\eta_\varepsilon*h)| d\mu(x)\\
	&=\frac{1}{\mu(B)}\int_B\Big|\int_{|z|<\varepsilon/8}{\eta_\varepsilon}(z)\Big[h(x-z) -(\Pi_{B-z}h)(x-z)\Big]\, dz\Big| d\mu(x)\\
	&\le \int_{|z|<\varepsilon/8}{\eta_\varepsilon}(z)\frac{1}{\mu(B)}\int_B\Big|\Big[h(x-z) -(\Pi_{B-z}h)(x-z)\Big]\Big| d\mu(x)\, dz\\
	&\le \int_{|z|<\varepsilon/8}{\eta_\varepsilon}(z)\frac{1}{\mu(B)}\int_{B-z}\Big|\Big[h(u) -(\Pi_{B-z}h)(u)\Big]\Big| e^{2z_1}d\mu(u)\, dz\\
	&\lesssim \|h\|_{{\BMO(\mu)}}\\
	\lesssim \beta.
\end{aligned}
\]
 
 Moreover, \eqref{eq:approx-local-remainder} and \eqref{eq:local-volume} imply
\[
 \sup_c\int_{B(c,1)}|h(x)|\,dx\lesssim \sup_c\frac{1}{\mu(B(c,1))}\int_{B(c,1)}|h(x)|\,d\mu(x)\le  C\beta.
\]
Thus $\|\eta_\varepsilon*h\|_\infty\le
C\|\eta_\varepsilon\|_\infty\beta<\infty$ for fixed
$\varepsilon>0$.

Finally, let $g_m=F_m+\eta_{1/m}*h$. These functions are bounded
and have uniformly bounded BMO seminorms. By
Theorem~\ref{thm:JN}, $h\in L^s_{\mathrm{loc}}$ for every finite
$s$; hence $\eta_{1/m}*h\to h$ in $L^s_{\mathrm{loc}}$.
The cutoffs give $F_m\to F$ locally, so
$g_m\to F+h=g-p$ in $L^s_{\mathrm{loc}}$.

{For $n=1$, the cutoff function $\chi_\perp$ is not needed.}

This completes our proof.
\end{proof}

\begin{theorem}\label{thm:duality}
The dual of ${H^1(\mu)}$ is ${\BMO(\mu)}/\cP$. Precisely, {every $g\in\BMO(\mu)$ defines a functional}
\[
 \Lambda_g(h)=\int_{\R^n}h(x)\overline{g(x)}\,d\mu(x)
\]
{initially on finite atomic sums. This functional extends uniquely to $H^1(\mu)$, and}
\begin{equation*}
 {\|\Lambda_g\|_{({H^1(\mu)})^*}\simeq\|g\|_{{\BMO(\mu)}}}.
\end{equation*}
Conversely, every bounded linear functional on ${H^1(\mu)}$ has this form. The
representing function is unique modulo $\cP$.
\end{theorem}

\begin{proof}
 Let $g\in{\BMO(\mu)}$, and take $g_m$ and $p$ from Lemma~\ref{lem:bounded-approximation}. Fix $b>0$ and $1<q<\infty$. If $a$ is a
$(1,q;\cP)_b$-atom supported in $B$, then H\"older's inequality,
Proposition~\ref{prop:BMO-scale} and Theorem~\ref{thm:JN} give
\[
 \left|\int a\overline{g_m}\,d\mu\right|
 =\left|\int a\overline{(g_m-\Pi_Bg_m)}\,d\mu\right|
 \le C_{b,q}\|g\|_{{\BMO(\mu)}}.
\]
For an arbitrary infinite atomic expansion
$h=\sum_k\lambda_ka_k$ in $L^1(\mu)$, since   $g_m\in L^\infty(\mu)$, the  integration is justified. Taking the infimum over all such expansions yields
\[
 \left|\int h\overline{g_m}\,d\mu\right|
 \le C_{b,q}\|g\|_{{\BMO(\mu)}}\|h\|_{{H^{1,q}_b(\mu)}}.
\]
If $h$ is a finite sum of $q$-atoms, it has bounded support and belongs
to $L^q$. By Lemma~\ref{lem:bounded-approximation} and $\int h\overline p\,d\mu=0$, we have
\begin{equation}\label{eq:bmo-pairing-infinite-norm}
 \left|\int h\overline g\,d\mu\right| = {\left|\int h\,\overline{(g-p)}\,d\mu\right|}=\lim_{m\to \infty}\left|\int h\overline g_m\,d\mu\right|
 \le C_{b,q}\|g\|_{{\BMO(\mu)}}\|h\|_{{H^{1,q}_b(\mu)}}.
\end{equation}
{This proves that} {the function $\Lambda_g$ defines a bounded linear functional} on ${H^1(\mu)}$. 
 
 \bigskip
 
Conversely, let $\Lambda\in({H^1(\mu)})^*$ and put $B_R=B(0,R)$ for
integers $R\ge1$. Let
\[
 V_R=\left\{u\in L^2(B_R,\mu):\int_{B_R}up\,d\mu=0\quad \ \ \text{for all} \ p\in\cP\right\}.
\]
We can cover $B_R$ by $N_R$ {overlapping balls} of a fixed radius at most $1$. The constant $C_R$ may depend on $R$ (through $N_R$ and the
ratio of radii), but not on $u$. Then we have
\[
 \|u\|_{{H^1(\mu)}}\le \mu(B_R)^{1/2} \|u\|_{L^2(B_R,\mu)}\le C_R\|u\|_{L^2(B_R,\mu)},\qquad u\in V_R.
\]
Since \(V_R\) is a closed subspace of \(L^2(B_R,\mu)\), the Riesz representation theorem gives \(f_R\in L^2(B_R,\mu)\) such that
\[
\Lambda(u)=\int_{B_R}u\overline{f_R}\,d\mu,
\qquad u\in V_R.
\]
Moreover, \(f_R\) is unique up to an element of \(\mathcal P\), since \(V_R^\perp=\mathcal P|_{B_R}\).

Now, if \(u\in V_R\), then its extension by zero belongs to \(V_{R+1}\). Hence
\[
\int_{B_R}u\overline{(f_{R+1}-f_R)}\,d\mu=0,
\qquad u\in V_R.
\]
It follows that \(f_{R+1}-f_R\in\mathcal P\) on \(B_R\). By subtracting a suitable element of \(\mathcal P\) from \(f_{R+1}\), we may assume that \(f_{R+1}=f_R\) on \(B_R\). Repeating this argument, we obtain a function \(f\in L^2_{\mathrm{loc}}(\mu)\) such that \(f=f_R\) on \(B_R\) for every \(R\). Consequently,
\[
\Lambda(u)=\int u\overline f\,d\mu
\]
for every compactly supported \(u\in L^2(\mu)\) with vanishing \(\mathcal P\)-moments, and hence for every finite linear combination of \((1,2;\mathcal P)_1\)-atoms.

Let $B$ be a ball of radius at most $1$, and put $g_B=f-\Pi_Bf$. If $g_B\ne0$, then
\[
 a_B=\frac{g_B\mathbf1_B}
              {\mu(B)^{1/2}\|g_B\|_{L^2(B,\mu)}}
\]
is an $(1,2;\cP)_1$-atom. Hence,
\[
 \mu(B)^{-1/2}\|f-\Pi_Bf\|_{L^2(B,\mu)}
 =|\Lambda(a_B)|\le\|\Lambda\|.
\]
Thus ${\|f\|_{{\BMO_1^2(\mu)}}}\le\|\Lambda\|$. The bounded functional
constructed from $f$ in the first part agrees with $\Lambda$ {on finite
linear combinations of $(1,2;\cP)_1$-atoms}, and hence on ${H^1(\mu)}$ by density. The integral expression  for every finite-$q$ atomic decomposition follows from the first part and the uniqueness modulo $\cP$ follows from zero local oscillation.

This completes our proof.
\end{proof}

We also define ${H^{1,q}_{b,\mathrm{fin}}(\mu)}$ as the space of {finite
linear combinations of $(1,q;\cP)_b$ atoms}, equipped with the norm
\[
\|f\|_{{H^{1,q}_{b,\mathrm{fin}}(\mu)}}
=\inf\left\{\sum_{k=1}^N|\lambda_k|:
f=\sum_{k=1}^N\lambda_ka_k,\ N<\infty\right\},
\]
{where the infimum is taken over all finite atomic decompositions.}

\begin{proposition}\label{prop:finite-atomic-norm}
Let $b>0$ and $1<q<\infty$. For every
$f\in {H^{1,q}_{b,\mathrm{fin}}(\mu)}$,
\begin{equation*}
 \|f\|_{{H^{1,q}_b(\mu)}}
 \le \|f\|_{{H^{1,q}_{b,\mathrm{fin}}(\mu)}}
 \le C_{b,q,n}\|f\|_{{H^{1,q}_b(\mu)}}.
\end{equation*}

\end{proposition}

\begin{proof}
The first inequality is {obvious}; hence, we need only to prove the second inequality. Let $\Lambda$ be a bounded
linear functional on ${H^{1,q}_{b,\mathrm{fin}}(\mu)}$. If $B$ is a ball with  radius at
most $b$ and $u\in L^q(B,\mu)$ annihilates $\cP$, then $u$ is a multiple
of an $(1,q;\cP)_b$-atom. Thus
\begin{equation}\label{eq:finite-functional-local}
 |\Lambda(u)|\le \|\Lambda\|\mu(B)^{1/q'}\|u\|_{L^q(B,\mu)}.
\end{equation}
For $R\ge1$, set $B_R=B(0,R)$. If $u\in L^q(B_R,\mu)$ annihilates $\cP$, {by a covering argument with balls of radius $b$} and
Lemma~\ref{lem:moment-transfer} as in Proposition \ref{prop:H1-independence}, we have
\[
 \|u\|_{{H^{1,q}_{b,\mathrm{fin}}(\mu)}}
 \le C_{R,b,q}\|u\|_{L^q(B_R,\mu)}.
\]
Thus $\Lambda$ is continuous on this closed subspace of $L^q(B_R,\mu)$. The duality of $L^q(B_R,\mu)$ with $L^{q'}(B_R,\mu)$ gives a representative
$g_R\in L^{q'}(B_R,\mu)$, unique modulo $\cP$. As in the proof of
Theorem~\ref{thm:duality}, there exists $g\in L^{q'}_{\mathrm{loc}}$ such that
\[
 \Lambda(f)=\int f\overline g\,d\mu
\]
for every $f\in {H^{1,q}_{b,\mathrm{fin}}(\mu)}$.

By \eqref{eq:finite-functional-local} and   $L^q(\mu)$--$L^{q'}(\mu)$ duality,
for every ball $B$ of radius at most $b$,
\[
 \inf_{p\in\cP}\|g-p\|_{L^{q'}(B,\mu)}
 \le \|\Lambda\|\mu(B)^{1/q'}.
\]
This, together with Lemma~\ref{lem:projection}, Proposition~\ref{prop:BMO-scale} and
Theorem~\ref{thm:JN}, gives $\|g\|_{{\BMO(\mu)}}\le C\|\Lambda\|$.
Similarly to \eqref{eq:bmo-pairing-infinite-norm}, we have, for every $f\in {H^{1,q}_{b,\mathrm{fin}}(\mu)}$,
\[
 |\Lambda(f)|
 =\left|\int f\overline g\,d\mu\right|
 \le C\|\Lambda\|\|f\|_{{H^{1,q}_b(\mu)}}.
\]
By the Hahn--Banach theorem, applied to the normed vector space
${H^{1,q}_{b,\mathrm{fin}}(\mu)}$,
\[
 \|f\|_{{H^{1,q}_{b,\mathrm{fin}}(\mu)}}
 {=\sup_{\|\Lambda\|\le1}|\Lambda(f)|\le C\|f\|_{{H^{1,q}_b(\mu)}}}.
\]

This completes our proof.
\end{proof}

\subsection{Interpolation Theorems}\label{sec:interpolation}

\paragraph{An auxiliary atomic space.}

Let $\mathcal D$ be the dyadic lattice obtained by subdividing the unit
cubes $k+[0,1)^n$, $k\in\mathbb Z^n$. An auxiliary $(1,q)$-atom is of one of
the following two types:
\begin{enumerate}[label=(\roman*)]
\item it is supported in a proper dyadic subcube $Q$ of a unit cube,
      satisfies $\int a\,d\mu=0$, and
      $\|a\|_{L^q(\mu)}\le\mu(Q)^{1/q-1}$;
\item it is supported in a unit dyadic cube $Q$, no cancellation is
      required, and $\|a\|_{L^q(\mu)}\le\mu(Q)^{1/q-1}$.
\end{enumerate}
We define the  space $X^{1}_q$ as a set of all functions $f\in L^1(\mu)$ such that
$$
	f=\sum_{k=1}^\infty\lambda_ka_k\quad\text{in }L^1(\mu),
	\qquad \sum_{k=1}^\infty|\lambda_k|<\infty,
$$
where each $a_k$ is a $(1,q)$-atom. We also define 
\[
\|f\|_{X^{1}_q}
=\inf\left\{\sum_{k=1}^\infty|\lambda_k|:
f=\sum_{k=1}^\infty\lambda_ka_k\right\}.
\]

\begin{lemma}\label{lem:aux-X1-q}
The spaces $X^1_q$, $1<q\le\infty$, coincide with equivalent norms.
{Hence, we define $X^1$ as any of the spaces $X^1_q$} for $1<q\le \infty$.
\end{lemma}

\begin{proof} It suffices to show that $X^1_q = X^1_\infty$ for all $1<q<\infty$. Fix $1<q<\infty$. The inclusion $X^{1}_\infty\hookrightarrow X^1_q$ is obvious by  H\"older's inequality. It remains to show that $X^{1}_q\hookrightarrow X^1_\infty$.  Let $a$ be an $(1,q)$-atom supported in a dyadic cube $Q$ with $\ell(Q)\le 1$. Choose $K>1$ and
denote by $\{Q_i \}_{i\in \mathcal I}$ the set of all maximal  dyadic subcubes $Q_i\subset Q$ satisfying 
\[
 \frac1{\mu(Q_i)}\int_{Q_i}|a|^q\,d\mu
 >\left(\frac K{\mu(Q)}\right)^q,
\]
which implies
\[
 \sum_i\mu(Q_i)\le K^{-q}\mu(Q).
\]
Since the parent of $Q_i$ does not satisfy the stopping time condition, we have
\[
 \frac1{\mu(Q_i)}\int_{Q_i}|a|^q\,d\mu
 \le C\left(\frac K{\mu(Q)}\right)^q.
\]
For each $i\in \mathcal I$, define
\[
 b_i=(a-a_{Q_i})\mathbf1_{Q_i},
 \qquad
 a_{Q_i}=\frac1{\mu(Q_i)}\int_{Q_i}a\,d\mu.
\]
{Set}
\[
g=a-\sum_i b_i.
\]
It is easy to see that 
\[
{\int b_i\,d\mu=0.}
\]
In addition,
\[
\begin{aligned}
	\sum_i\mu(Q_i)^{1-1/q}\|b_i\|_{L^q(\mu)}&\le \Big(\sum_i\mu(Q_i)\Big)^{1/q'}\Big(\sum_i \|b_i\|^q_{L^q(\mu)}\Big)^{1/q}\\
	&\le 2K^{-(q-1)}{\mu(Q)^{1/q'} \|a\|_{L^q(\mu)}}\\
	&\le  2K^{-(q-1)}.
\end{aligned}
\]

By the dyadic Lebesgue differentiation theorem, we have
\[
\|g\|_\infty\le CK/\mu(Q).
\]
This, in combination with the fact  $\supp g\subset Q$, implies that $g/(CK)$
is an $(1,\infty)$ atom of type (i) if $\ell(Q)<1$,
and of type (ii) if $\ell(Q)=1$. Hence,
\[
\|g\|_{X^1_\infty}\le CK.
\]

We have thus proved that every $(1,q)$ atom $a$
associated with a cube $Q$ can be written as
\[
a=g+\sum_i\lambda_i a_i,
\]
where $g\in X^1_\infty$ satisfies
$\|g\|_{X^1_\infty}\le CK$, each $a_i$ is an
$(1,q)$ atom of type (i), and
\[
\sum_i|\lambda_i|\le\frac12.
\]

At this stage, arguing similarly to the proof of Proposition \ref{prop:H1-independence}, we complete our proof.
\end{proof}

We need an operator which converts auxiliary $(1,q)$ atoms into atoms annihilating
$\cP$, without losing the ability to recover a function in
$L^p(\mu)$, $p>1$. A translation in the drift direction supplies
both properties: after choosing its weight to preserve the $\mu$-density,
{the difference between the identity and a forward average has zero
moments against $\cP$.} At the same time the average is a strict
contraction on $L^p(\mu)$ because the measure grows in that direction.

Choose a nonnegative radial function $\zeta\in C_c^\infty(B(0,1/8))$ with
$\displaystyle \int \zeta(z)\,dz=1$. Define
\begin{equation}\label{eq:Aop}
 \mathcal Af(x)=\int \zeta(z)e^{-2(\frac12+z_1)}
 f(x-\frac12e_1-z)\,dz,
 \qquad \text{and} \ \ 
 \mathcal S=I-\mathcal A.
\end{equation}
If $F_f(x)=e^{2x_1}f(x)$, then
\begin{equation}\label{eq:A-density}
 e^{2x_1}\mathcal Af(x)
 =(\zeta*F_f)(x-\frac12e_1).
\end{equation}
Since $\displaystyle  \int z_j\zeta(z)\,dz=0$ for $j\ge2$, using  Fubini's theorem, we have
\begin{equation}\label{eq:A-moment-preserve}
 \int \mathcal Af\,d\mu=\int f\,d\mu,
 \qquad
 \int x_j\mathcal Af\,d\mu=\int x_jf\,d\mu,
 \quad j\ge2.
\end{equation}
{Consequently,}
\[
\int \mathcal S f(x) p(x)\,d\mu(x)=0 \ \ \ \ \ \text{for all $p\in \cP$}.
\]

\begin{lemma}\label{lem:S-properties}
The operator $\mathcal S$ is bounded on $L^1(\mu)$ and {its inverse is bounded on $L^p(\mu)$} for every $1<p<\infty$. In addition,
\begin{equation*}
 \mathcal S:X^1\longrightarrow {H^1(\mu)}
\end{equation*}
is bounded.
\end{lemma}

\begin{proof}
For $h_1>0$,
\[
 \|e^{-2h_1}f(\cdot-h)\|_{L^p(\mu)}
 =e^{-2h_1(1-1/p)}\|f\|_{L^p(\mu)}.
\]
On the support of $\zeta$,  we have $\frac12+z_1\ge\frac38$. By Minkowski's inequality, we have
\begin{equation*}
 \|\mathcal A\|_{L^p(\mu)\to L^p(\mu)}
 \le e^{-\frac34(1-1/p)}<1,
 \qquad 1<p<\infty.
\end{equation*}
For $p=1$, $ \|\mathcal A\|_{L^p(\mu)\to L^p(\mu)}\le 1$, hence $\mathcal S$ is bounded on $L^1(\mu)$. For $1<p<\infty$, we have $\|\mathcal A\|_{L^p(\mu)\to L^p(\mu)}<1$. Therefore, for $1<p<\infty$,
\[
 \mathcal S^{-1}=\sum_{m=0}^\infty\mathcal A^m
\]
converges in operator norm on $L^p(\mu)$.

Since $\mathcal S$ is bounded on $L^1(\mu)$, by Lemma~\ref{lem:aux-X1-q}, it suffices to prove that
\begin{equation}\label{eq-Sa less 1}
\|\mathcal S a\|_{{H^1(\mu)}}\le C
\end{equation}
for all $(1,2)$-atoms.

{Let $a$ be an $(1,2)$-atom supported in a dyadic cube $Q$ centered at $c$} and $\ell(Q)\le 1$. The map
$\tau_c f(x)=e^{2c_1}f(x+c)$ preserves the $L^1(\mu)$ norm, commutes with $\mathcal A$. Therefore, we might assume that $c=0$. If $\ell(Q)=1$, i.e.,
$a$ is an $(1,2)$ atom of  type (ii).   By
\eqref{eq:A-moment-preserve}, 
\[
\int \mathcal Sa(x)p(x)\,d\mu(x)={\int a(x)p(x)\,d\mu(x)-\int\mathcal Aa(x)p(x)\,d\mu(x)} =0  \ \ \ \ \text{for all $p\in \cP$}.
\]
Moreover,  from \eqref{eq:Aop}, we have $\supp \mathcal Sa \subset  B(0,2\sqrt n)$. From the boundedness of $\mathcal S$ on $L^2(\mu)$, 
\[
\|\mathcal Sa\|_{L^2(\mu)}\lesssim \|a\|_{L^2(\mu)}\lesssim \mu(B(0,2\sqrt n))^{-1/2}.
\]
Hence, $\mathcal Sa$ is a multiple of an $(1,2;\cP)_{2\sqrt{n}}$ atom. This, together with Proposition~\ref{prop:H1-independence}, implies that 
$\|\mathcal Sa\|_{{H^1(\mu)}}\le C$.

Let  $a$ be an $(1,2)$ atom of  type (i), i.e.,  $r:=\ell(Q)<1$. If
$r\ge1/8$, {arguing as in the case $r=1$}, we also obtain \eqref{eq-Sa less 1}. It remains to prove the inequality when $r<1/8$. Set
\[
 F_a=e^{2x_1}a(x),
 \qquad
 M_j=\int x_ja(x)\,d\mu(x),\quad 2\le j\le n.
\]
Then $\displaystyle \int F_a\,dx=0$, $\|F_a\|_{L^1(dx)}\le1$, and 
\begin{equation}\label{eq-M}
	{|M|:=\sum_{j=2}^n|M_j|\le Cr.} 
\end{equation}
From
\eqref{eq:A-density} and the mean-zero condition,
\[
 |(\zeta*F_a)(u)|
 =\left|{\int_{Q}[\zeta(u-y)-\zeta(u)]F_a(y)\,dy}\right|
 \le Cr.
\]

From  \eqref{eq:A-density} and the facts ${\supp F_a\subset B(0,\sqrt n\,r)}$ and $\supp\zeta\subset B(0,1/8)$, we have ${\supp \mathcal Aa\subset B(0,2\sqrt n)}$. This, together with \eqref{eq:A-density} and the {above inequality}, implies 
\begin{equation}\label{eq:A-small-L2}
 \|\mathcal Aa\|_{L^2(\mu)}\le  Cr.
\end{equation}

Choose $\varphi\in C_c^\infty(B(0,1))$ with $\displaystyle \int\varphi\,dx=1$ and set
$\varphi_\rho(x)=\rho^{-n}\varphi(x/\rho)$. For $r\le\rho\le1$, define
\begin{equation*}
 h_\rho(x)=-e^{-2x_1}\sum_{j=2}^nM_j\partial_j\varphi_\rho(x).
\end{equation*}
Then, by integration by parts, 
\begin{equation}\label{eq:h-rho-moments}
 \int h_\rho\,d\mu=0,
 \qquad
 \int x_jh_\rho\,d\mu=M_j,  \ \ \ 2\le j\le n.
\end{equation}
Moreover,
\begin{equation}\label{eq:h-rho-size}
 \|h_\rho\|_{L^2(\mu)}
 \le C\frac{|M|}{\rho}\mu(B(0,\rho))^{-1/2}.
\end{equation}
Let $N\in \mathbb N$ be the smallest number satisfying  $2^Nr\ge1/2$ and set $\rho_k=2^kr$ for $k=0,\ldots,N$. Then, we write
\begin{equation}\label{eq:S-telescope}
 \mathcal Sa
 =(a-h_r)+\sum_{k=0}^{N-1}(h_{\rho_k}-h_{\rho_{k+1}})
   +(h_{\rho_N}-\mathcal Aa).
\end{equation}
Each term on the right hand side of \eqref{eq:S-telescope} annihilates  $\cP$ due to
\eqref{eq:A-moment-preserve} and \eqref{eq:h-rho-moments}. For the first term $a-h_r$, we have ${\supp(a-h_r)\subset B(0,\sqrt n\,r)}$; moreover, by \eqref{eq:h-rho-size} and \eqref{eq-M},
\[
\|a-h_r\|_{L^2(\mu)}\le \|a\|_{L^2(\mu)}+\|h_r\|_{L^2(\mu)}\lesssim {(1+\frac{|M|}{r})} \mu(B(0,r))^{-1/2}\lesssim \mu(B(0,r))^{-1/2}.
\]
Therefore, $a-h_r$ is a multiple of an $(1,2)$ atom, which implies $\|a-h_r\|_{X^1}\lesssim 1$. 

For each $k=0,\ldots, N-1$, we have $\supp(h_{\rho_k}-h_{\rho_{k+1}})\subset B(0,2^{k+1}r)$. In addition, by \eqref{eq:h-rho-size}, \eqref{eq-M} and \eqref{eq-local doubling},
\[
\|h_{\rho_k}-h_{\rho_{k+1}}\|_{L^2(\mu)}{\lesssim \frac{|M|}{2^kr}}\mu(B(0,2^{k+1}r))^{-1/2}\lesssim 2^{-k}\mu(B(0,2^{k+1}r))^{-1/2}.
\]
It follows that 
\[
\|h_{\rho_k}-h_{\rho_{k+1}}\|_{X^1}\lesssim 2^{-k},
\]
which implies
\[
\Big\|\sum_{k=0}^{N-1}(h_{\rho_k}-h_{\rho_{k+1}})\Big\|_{X^1}\lesssim 1.
\]

{The last term $h_{\rho_N}- \mathcal A a$ is supported in $B(0,2\sqrt n)$. Moreover, by} \eqref{eq:A-small-L2}, \eqref{eq:h-rho-size}, \eqref{eq-M} and $r<1/8$, we have
\[
{\|h_{\rho_N}-\mathcal Aa\|_{L^2(\mu)}}\lesssim {\frac{|M|}{\rho_N} + r \lesssim 1,}
\]
which implies
\[
\|h_{\rho_N}- \mathcal A a\|_{X^1}\lesssim 1.
\]
This completes our proof.
\end{proof}

\begin{lemma}\label{lem:dyadic-interpolation}
Let $T$ be a sublinear operator, bounded from $X^1$ to $L^1(\mu)$
and on $L^2(\mu)$. Then $T$ is bounded on $L^p(\mu)$ for every
$1<p<2$.
\end{lemma}

\begin{proof}
Fix $1<q<p<2$. Let $\mathcal M_d$ be the dyadic maximal operator defined by
\[
\mathcal M_d f(x) ={\sup_{\substack{Q\ni x\\ Q\in \mathcal D}}}\frac{1}{\mu(Q)}\int_Q |f(y)|d\mu(y), \ \ \ \ x\in \mathbb R^n,
\]
where the supremum is taken over all dyadic cubes $Q$ such that $x\in Q$ and $\ell(Q)\le 1$.

Due to \eqref{eq-local doubling}, $\mathcal M_d$ is bounded on $L^r(\mu)$ for $1<r<\infty$. Setting  $\mathcal M_{d,q}f=(\mathcal M_d(|f|^q))^{1/q}$, then
\[
 \|\mathcal M_{d,q}f\|_{L^p(\mu)}\le C_{p,q}\|f\|_{L^p(\mu)}.
\]
Let $f\in L^\infty(\mu)$ with  compact  support.
{ 
For $\lambda>0$, let
\[
 \mathcal D_\lambda
 =\left\{Q\in\mathcal D:
 \frac1{\mu(Q)}\int_Q|f|^q\,d\mu>\lambda^q\right\}.
\]
Let
$\mathcal Q_\lambda$ be the family of maximal cubes in
$\mathcal D_\lambda$, and set
$\Omega_\lambda=\bigcup_{Q\in\mathcal Q_\lambda}Q$.
If $\mathcal D_\lambda$ is empty, then $\mathcal Q_\lambda$ is also empty,
$\Omega_\lambda=\emptyset$.
The
maximal cubes are pairwise disjoint. Obviously,
\[
 \Omega_\lambda=\{x:\mathcal M_{d,q}f(x)>\lambda\}.
\]
}
Set
\[
 b_\lambda=\sum_{Q\in\mathcal Q_\lambda}b_Q,
 \qquad g_\lambda=f-b_\lambda,
 \qquad
 b_Q=\begin{cases}
 (f-f_Q)\mathbf1_Q,&\ell(Q)<1,\\
 f\mathbf1_Q,&\ell(Q)=1,
 \end{cases}
\]
where $f_Q=\mu(Q)^{-1}\int_Qf\,d\mu$. If $\mathcal Q_\lambda=\emptyset$, then we set $b_\lambda =0$ and $g_\lambda = f$.

{ 
If $\ell(Q)<1$, its dyadic parent $\widehat Q$
belongs to $\mathcal D$ and does not belong to $\mathcal D_\lambda$.
By local doubling property \eqref{eq-local doubling},
\[
\begin{aligned}
	\|b_Q\|_{L^q(\mu)}
	&\le C\left(\int_Q|f|^q\,d\mu\right)^{1/q}\\
	&\le C\mu(\widehat Q)^{1/q}\left(\frac{1}{\mu(\widehat Q)}\int_{\widehat Q}|f|^q\,d\mu\right)^{1/q}\\
	& \le C\lambda\mu( Q)^{1/q}\\
	&\le C\mu( Q)^{1/q} \inf_{x\in Q} \mathcal M_{d,q}f(x).
\end{aligned}
\]
Moreover, $\int_Q b_Q\,d\mu=0$. {Hence, $b_Q$ is a multiple of an auxiliary $(1,q)$-atom when $\ell(Q)<1$.}

If $\ell(Q) =1$, then we have
\[
\begin{aligned}
	\|b_Q\|_{L^q(\mu)}&=  \mu(Q)^{1/q}\left(\frac{1}{\mu(Q)}\int_{{Q}}|f|^q\,d\mu\right)^{1/q}\\
	 &\le \mu( Q)^{1/q} \inf_{x\in Q} \mathcal M_{d,q}f(x),
\end{aligned}
\]
{which implies that $b_Q$ is a multiple of an auxiliary $(1,q)$-atom when $\ell(Q)=1$.}

Therefore,
\[
\begin{aligned}
 \|b_Q\|_{X^1_q}
 &\le \mu(Q)^{1-1/q}\|b_Q\|_{L^q(\mu)}\\
 &\le C\mu(Q)  \inf_{x\in Q} \mathcal M_{d,q}f(x)\\
& \le C\int_Q  \mathcal M_{d,q}f\,d\mu.
\end{aligned}
\]
Since the family $\mathcal Q_\lambda$ is pairwise disjoint and
$\mathcal M_{d,q}f>\lambda$ on each cube in $\mathcal Q_\lambda$, Lemma~\ref{lem:aux-X1-q} gives
\begin{equation*}
 \|b_\lambda\|_{X^1}
 \le C\sum_{Q\in\mathcal Q_\lambda}\int_Q\mathcal M_{d,q}f\,d\mu
 \le C\int_{\{\mathcal M_{d,q}f>\lambda\}}\mathcal M_{d,q}f\,d\mu.
\end{equation*}
}
In addition, we have $g_\lambda=f$ on $\mathbb R^n\setminus \Omega_\lambda$, while $g_\lambda=f_Q$ when $Q\in \mathcal{Q}_\lambda, \ell(Q)<1$,    and  $g_\lambda=0$ on $Q\in \mathcal{Q}_\lambda, \ell(Q)=1$.  Consequently,
\begin{equation*}%
	\|g_\lambda\|_{L^2(\mu)}^2
	\le C\Big(\int_{\{\mathcal M_{d,q}f\le\lambda\}}[\mathcal M_{d,q}f]^2\,d\mu
	+ {\sum_{\substack{Q\in\mathcal Q_\lambda\\ \ell(Q)<1}}}\int_{Q} |f_Q|^2d\mu \Big).
\end{equation*}
For $\ell(Q)<1$, its dyadic parent $\widehat Q$
belongs to $\mathcal D$ and does not belong to $\mathcal D_\lambda$.
By local doubling property \eqref{eq-local doubling},
\[
\begin{aligned}
	|f_Q|
	&\le  \Big(\frac{\mu(\widehat Q)}{\mu(Q)}\Big)^{1/q}\Big(\frac{1}{\mu(\widehat Q)}\int_{\widehat Q}|f|^q\,d\mu\Big)^{1/q} \le C\lambda.
\end{aligned}
\]
{These inequalities imply} 
\begin{equation*}
 \|g_\lambda\|_{L^2(\mu)}^2
 \le C\left(\int_{\{\mathcal M_{d,q}f\le\lambda\}}[\mathcal M_{d,q}f]^2\,d\mu
   +\lambda^2\mu\{\mathcal M_{d,q}f>\lambda\}\right).
\end{equation*}
{The $X^1$--$L^1(\mu)$ boundedness}, {the $L^2(\mu)$-boundedness}  and sublinearity of $T$ yield
\[
\begin{aligned}
\mu\{| {T}f|>2\lambda\}
&\le \frac{C}{\lambda}\|b_\lambda\|_{X^1}
  +\frac{C}{\lambda^2}\|g_\lambda\|_{L^2(\mu)}^2\\
&\le C\left(\frac1\lambda\int_{\{\mathcal M_{d,q}f>\lambda\}}\mathcal M_{d,q}f\,d\mu
 +\frac1{\lambda^2}\int_{\{\mathcal M_{d,q}f\le\lambda\}}[\mathcal M_{d,q}f]^2\,d\mu
 +\mu\{\mathcal M_{d,q}f>\lambda\}\right).
\end{aligned}
\]
Multiplying  by $\lambda^{p-1}$, integrating over $(0,\infty)$ and using 
Fubini's theorem, we obtain
\[
 \|{T}f\|_{L^p(\mu)}^p\le C_{p,q}{\|\mathcal M_{d,q}f\|_{L^p(\mu)}^p}
 \le C_{p,q}\|f\|_{L^p(\mu)}^p.
\]
This completes the proof.
\end{proof}

\begin{theorem}\label{thm:Hardy-interpolation}
Let $T$ be sublinear and bounded from ${H^1(\mu)}$ to $L^1(\mu)$ and
on $L^2(\mu)$. Then $T$ is bounded on $L^p(\mu)$ for every
$1<p<2$.
\end{theorem}

\begin{proof}
Let $\mathcal S$ be defined by \eqref{eq:Aop}. By
Lemma~\ref{lem:S-properties}, $T\mathcal S$ is bounded from $X^1$
to $L^1(\mu)$ and on $L^2(\mu)$. Lemma~\ref{lem:dyadic-interpolation} gives
\[
 \|T\mathcal Sh\|_{L^p(\mu)}\le C_p\|h\|_{L^p(\mu)},\qquad 1<p<2.
\]
Since $\mathcal S^{-1}$ is bounded on $L^p(\mu)$,
\[
 \|Tf\|_{L^p(\mu)}\le C_p\|\mathcal S^{-1}f\|_{L^p(\mu)}
 \le C_p\|f\|_{L^p(\mu)}.
\]
This completes our proof.
\end{proof}

As a direct consequence of Theorem \ref{thm:Hardy-interpolation} and Theorem \ref{thm:duality}, we have:
\begin{theorem}\label{thm:BMO-interpolation}
Let $T$ be linear and bounded on $L^2(\mu)$. Assume that
$T:L^\infty(\mu)\to{\BMO(\mu)}/\cP$ is bounded, with the endpoint
action agreeing with $T$ on $L^2(\mu)\cap L^\infty(\mu)$.
Then $T$ is bounded on $L^p(\mu)$ for every $2<p<\infty$.
\end{theorem}

\section{Riesz transforms on Hardy and BMO spaces}\label{sec:Riesz}
Set $Uf=e^{x_1}f$. Then $U:L^2(\mu)\to L^2(dx)$ is unitary and
\begin{equation}\label{eq:ground-state}
	U\cL U^{-1}=-\Delta+1.
\end{equation}
Hence, the heat kernel of $e^{-t\cL}$ is defined by
\begin{equation}\label{eq:heat-density}
 e^{2x_1}p_t(x,y)=Q_t(x-y),
 \qquad
 Q_t(z)=(4\pi t)^{-n/2}
 \exp\left(-\frac{|z-2te_1|^2}{4t}\right).
\end{equation}
Set
 \begin{equation}\label{eq:Phi}
 \Phi(z)=\frac1{\sqrt\pi}\int_0^\infty t^{-1/2}Q_t(z)\,dt.
\end{equation}
Set $\Phi(z)=e^{z_1}G(z)$. Then
\begin{equation}\label{eq:G}
 G(z)=d_n\int_0^\infty t^{-(n+1)/2}
 e^{-t-|z|^2/(4t)}\,dt,
 \qquad d_n=\frac1{\sqrt\pi(4\pi)^{n/2}}.
\end{equation}

Denote by $K_j$ and $K_j^*$  the kernels of $R_j$ and $R_j^*$ relative to
$d\mu$. Then, we have
\begin{equation}\label{eq:Fj}
 e^{2x_1}K_j(x,y)=F_j(x-y),
 \qquad
 F_j=(\partial_j-2\delta_{1j})\Phi,
\end{equation}
and
\begin{equation}\label{eq:Fj-star}
 e^{2x_1}K_j^*(x,y)=F_j^*(x-y),
 \qquad
 F_j^*=-\partial_j\Phi.
\end{equation}
Indeed \eqref{eq:Fj} follows by differentiating
$e^{-2x_1}\Phi(x-y)$, and \eqref{eq:Fj-star} follows from
$R_j^*=-R_j-2\delta_{1j}\cL^{-1/2}$.

The Riesz transforms are well-defined for $C_c^\infty(\R^n)$ functions and bounded on $L^2(\mu)$.

\begin{lemma}\label{lem:test-Riesz-bounded}
If $\varphi\in C_c^\infty(\R^n)$, then
\[
 R_j\varphi,\ R_j^*\varphi\in L^\infty(\mu)\cap L^2(\mu),
 \qquad 1\le j\le n.
\]
\end{lemma}

\begin{proof}
{
By \eqref{eq:L2-Riesz}, both $R_j\varphi$ and $R_j^*\varphi$
belong to $L^2(\mu)$. It remains to prove that $ R_j\varphi,\ R_j^*\varphi\in L^\infty(\mu)$.

{Recall the function $G$ defined in \eqref{eq:G}. With $r=|z|$, since}
\[
 t+\frac{r^2}{4t}=r+\frac{(r-2t)^2}{4t},
\]
we may write $G(z)=d_ne^{-r}J_0(r)$, where
\[
 J_m(r)=\int_0^\infty t^{-(n+1)/2-m}
 \exp\left(-\frac{(r-2t)^2}{4t}\right)\,dt,
 \qquad m=0,1.
\]
Splitting the integral into $0<t<r/4$, $r/4\le t\le r$ and
$t>r$, we obtain, for $r\ge1$,
\[
 J_m(r)\lesssim r^{-n/2-m},\qquad m=0,1.
\]

It follows that
\[
 |G(z)|\lesssim e^{-r}r^{-n/2}.
\]
Differentiating under the integral sign gives
\[
 \partial_jG(z)=-\frac{d_nz_j}{2}
 \int_0^\infty t^{-(n+3)/2}e^{-t-r^2/(4t)}\,dt.
\]
Thus $|\nabla G(z)|\lesssim re^{-r}J_1(r)
\lesssim e^{-r}r^{-n/2}$. Consequently,
\begin{equation}\label{eq:Bessel-tail}
 |\partial^\alpha G(z)|\le C_\alpha e^{-|z|}|z|^{-n/2},
 \qquad |z|\ge1,\quad |\alpha|\le1.
\end{equation}
This, together with the fact
\[
 F_j(z)=e^{z_1}(\partial_jG(z)-\delta_{1j}G(z)),\qquad
 F_j^*(z)=-e^{z_1}(\partial_jG(z)+\delta_{1j}G(z)),
\]
implies
\begin{equation}\label{eq:F-tail-pointwise}
 |F_j(z)|+|F_j^*(z)|
 \le C e^{z_1-|z|}|z|^{-n/2},\qquad |z|\ge1.
\end{equation}

Let $K=\supp\varphi$. For $x\in \mathbb R^n$ with $\dist(x,K)\ge1$, using \eqref{eq:Fj}, \eqref{eq:Fj-star} and \eqref{eq:F-tail-pointwise},
\[
\begin{aligned}
 |R_j\varphi(x)|+|R_j^*\varphi(x)|
 &\le C\int_K e^{-x_1+y_1-|x-y|}|x-y|^{-n/2}
               |\varphi(y)|\,dy\\
 &\le C\int_K|\varphi(y)|\,dy,
\end{aligned}
\]
because $|x-y|\ge|x_1-y_1|$ and $|x-y|\ge1$.

On the compact set $\{x:\operatorname{dist}(x,K)\le 1\}$, recall that $Uf=e^{x_1}f$. Then
{$R_j=U^{-1}TU$ and $R_j^*=U^{-1}T^*U$}, where $T=(\partial_j-\delta_{1j})(-\Delta+1)^{-1/2}$ has symbol
$m_j(\xi)=(i\xi_j-\delta_{1j})(1+|\xi|^2)^{-1/2}$ with $|m_j|\le 1$. {Since $U\varphi\in C_c^\infty$,}
Fourier inversion gives
\[
\|T(U\varphi)\|_\infty+{{\|T^*(U\varphi)\|_\infty}}
\le (2\pi)^{-n}\bigl(\|\widehat{U\varphi}\|_1+{\|\widehat{U\varphi}\|_1}\bigr)<\infty ,
\]
and $e^{\pm x_1}$ is bounded on the compact set, so $R_j\varphi$ and $R_j^*\varphi$ are bounded .

This completes our proof.
}
\end{proof}

\begin{lemma}\label{lem:kernel-derivatives}
There is $C>0$ such that, for $0<|z|\le4$,
\begin{equation}\label{eq:local-derivative}
 |\nabla F_j(z)|+|\nabla F_j^*(z)|
 \le C|z|^{-n-1},
 \qquad 1\le j\le n.
\end{equation}
Moreover,
\begin{equation}\label{eq:global-integrability}
 \partial_1F_1,\quad
 \partial_k\partial_\ell F_1\ (k,\ell\ge2),\quad
 \nabla F_j\ (j\ge2),\quad
 \nabla F_j^*\ (1\le j\le n)
\end{equation}
belong to $L^1(\{|x|>1\},dx)$.
\end{lemma}

\begin{proof}
By \eqref{eq:heat-density}, for every multi-index $\alpha$,
\begin{equation}\label{eq-Qt}
 |\partial^\alpha Q_t(z)|
 \le C_\alpha t^{-(n+|\alpha|)/2}
       e^{-|z-2te_1|^2/(8t)},
 \qquad
 \|\partial^\alpha Q_t\|_{L^1(dx)}\le C_\alpha t^{-|\alpha|/2}.
\end{equation}
Since
\[
 \frac{|z-2te_1|^2}{8t}
   =\frac{|z|^2}{8t}-\frac{z_1}{2}+\frac t2,
\]
we have, {for $|z|\le 4$,}
\[
 |\partial^\alpha Q_t(z)|
\le C_\alpha t^{-(n+|\alpha|)/2}
e^{-|z-2te_1|^2/(8t)}\simeq C_\alpha t^{-(n+|\alpha|)/2}
e^{-|z|^2/(8t)}e^{-t/2}. 
\]
This, together with \eqref{eq:Phi}, implies, for $|\alpha|=1,2$ and $0<|z|\le4$,
\[
 |\partial^\alpha\Phi(z)|
 \le C\int_0^\infty t^{-(n+|\alpha|+1)/2}
          e^{-|z|^2/(8t)-t/2}\,dt
 \le C|z|^{1-n-|\alpha|}.
\]
This proves \eqref{eq:local-derivative}.

In order to prove \eqref{eq:global-integrability}, we first have
\[
\begin{aligned}
	|\partial^\alpha \Phi(z)|&\lesssim \int_0^{1/8}  t^{-1/2} |{\partial^\alpha Q_t(z)}|dt +\int_{1/8}^\infty  t^{-1/2} |{\partial^\alpha Q_t(z)}|dt,
\end{aligned}
\]
which implies
\[
\begin{aligned}
	\int_{|z|>1}|\partial^\alpha \Phi(z)|\,dz&\lesssim 	\int_{|z|>1}\int_0^{1/8}  t^{-1/2} |{\partial^\alpha Q_t(z)}|dtdz +\int_{1/8}^\infty  t^{-1/2} \|{\partial^\alpha Q_t(z)}\|_{L^1(dz)}dt\\
	&=: E_1 +E_2.
\end{aligned}
\]
By \eqref{eq-Qt}, we have $E_2\lesssim 1$ for all $|\alpha|\ge 2$.

{For the term $E_1$, note that for $0<t\le1/8$}, $|z|>1$, we have $|z-2te_1|\gtrsim |z|$. Hence, from \eqref{eq-Qt},
\[
|\partial^\alpha Q_t(z)|\le C_\alpha t^{-(n+|\alpha|)/2}
e^{-|z|^2/(ct)}.
\]
{Therefore, for $|z|>1$,}
\[
\begin{aligned}
	E_1&{\lesssim \int_0^{1/8}t^{-(n+|\alpha|+1)/2}
       \int_{|z|>1}e^{-|z|^2/(ct)}\,dz\,dt}\\
	&{\lesssim \int_0^{1/8}t^{-(|\alpha|+1)/2}e^{-c'/t}\,dt
       \lesssim 1.}
\end{aligned}
\]
This shows that  all the derivatives in \eqref{eq:global-integrability}, {except
 $\partial_1F_1$, belong} to $L^1(\{|z|>1\}, dx)$.

Finally, the heat equation for \eqref{eq:heat-density} gives
\[
 \partial_tQ_t=\Delta Q_t-2\partial_1Q_t.
\]
Integrating by parts in \eqref{eq:Phi}, for $z\ne0$ we obtain
\[
 \partial_1F_1(z)
 =\frac1{2\sqrt\pi}\int_0^\infty t^{-1/2}\partial_t Q_t(z)\,dt
  -\frac1{\sqrt\pi}\sum_{k=2}^n
       \int_0^\infty t^{-1/2}\partial_k^2Q_t(z)\,dt.
\]
Both terms belong to $L^1(\{|z|>1\},dx)$ by the same argument.

This completes our proof.
\end{proof}

\begin{theorem}\label{thm:H1-endpoint}
For every $1\le j\le n$,
\begin{equation*}
 \|R_jf\|_{L^1(\mu)}+\|R_j^*f\|_{L^1(\mu)}
 \le C\|f\|_{{H^1(\mu)}}.
\end{equation*}
Hence, for $1<p<\infty$,
\[
\|R_jf\|_{L^p(\mu)}+\|R_j^*f\|_{L^p(\mu)}
\le C_p\|f\|_{L^p(\mu)}.
\]
\end{theorem}

\begin{proof}
	{For $1<p<2$, the second part follows from the $H^1(\mu)$--$L^1(\mu)$ boundedness and Theorem~\ref{thm:Hardy-interpolation}. The case $2<p<\infty$ follows by duality. Hence, it suffices to prove that}
\[
 \|R_ja\|_{L^1(\mu)}+\|R_j^*a\|_{L^1(\mu)}
\le C
\]
for every $(1,2;\cP)_1$ atom associated with a ball $B=B(c,r)$ with $r\le 1$.

To do this, we write
 \[
  \|R_ja\|_{L^1(\mu)}+\|R_j^*a\|_{L^1(\mu)}\le  \|R_ja\|_{L^1(2B,\mu)}+\|R_j^*a\|_{L^1(2B,\mu)} +\|R_ja\|_{L^1(\mathbb R^n\setminus 2B,\mu)}+\|R_j^*a\|_{L^1(\mathbb R^n\setminus 2B,\mu)}
 \]

{By H\"older's inequality}, the $L^2(\mu)$-boundedness of $R_j$ and $R_j^*$ and the local doubling \eqref{eq-local doubling},
\[
\begin{aligned}
	\|R_ja\|_{L^1(2B,\mu)}+\|R_j^*a\|_{L^1(2B,\mu)}&\le \mu(2B)^{1/2}\Big[\|R_ja\|_{L^2(2B,\mu)}+\|R_j^*a\|_{L^2(2B,\mu)}\Big]\\
	&\lesssim \mu(B)^{1/2}\|a\|_{L^2(\mu)}\\
	&\lesssim 1.
\end{aligned}
\]

We now estimate the integrals outside $2B$. Set $E=\mathbb R^n\setminus2B$.
For $x\in E$, using \eqref{eq:Fj}, we have
\[
 R_ja(x) =e^{-2x_1}\int_B F_j(x-y)a(y)\,d\mu(y).
\]
Therefore,
\begin{equation}\label{eq:weighted-L1-density}
\begin{aligned}
 \int_E|R_ja(x)|\,d\mu(x)
 &=\int_E\left|\int_B F_j(x-y)a(y)\,d\mu(y)\right|dx\\
 &=\int_{E\cap B(c,3)}\!\left|\int_B F_j(x-y)a(y)\,d\mu(y)\right|dx\\
 &\quad+\int_{E\setminus B(c,3)}\!\left|\int_B F_j(x-y)a(y)\,d\mu(y)\right|dx
 =:H_1+H_2.
\end{aligned}
\end{equation}
Since $\displaystyle \int_B a\,d\mu=0$, we have
\[
H_1= \int_{E\cap B(c,3)}\!\left|\int_B \Big[F_j(x-y)-F_j(x-c)\Big]a(y)\,d\mu(y)\right|dx
\]
If $x\in E\cap B(c,3)$ and $y\in B$,
then $|x-c|>2r$ and
$|x-c-t(y-c)|\ge |x-c|/2$ for $0\le t\le1$.
By the mean value theorem and \eqref{eq:local-derivative},
\[
\begin{aligned}
 H_1
 &\le C\int_B |a(y)|\,|y-c|\,d\mu(y)
     \int_{2r<|x-c|<3}|x-c|^{-n-1}\,dx\\
 &\le C r\|a\|_{L^1(\mu)}\int_{2r}^3 s^{-2}\,ds
 \le C.
\end{aligned}
\]
The same estimate holds for $R_j^*a$ by \eqref{eq:Fj-star}.

For $|x-c|\ge3$, if $j\ge 2$, then {using the cancellation of $a$}, 
\[
\begin{aligned}
	H_2
	&\le \int_B |a(y)|\int_{|x-c|\ge3}
	|F_j(x-y)-F_j(x-c)|\,dx\,d\mu(y).
\end{aligned}
\]
Note that all points on the line segment from $x-y$ to $x-c$
have norm greater than $2$, and $\nabla F_j\in L^1(\{|z|>1\},dz)$. This, together with the mean value theorem and
Lemma~\ref{lem:kernel-derivatives}, implies
\[
\begin{aligned}
 H_2
 &\le C\int_B|a(y)|\,|y-c|\,d\mu(y)\\
& \le Cr\|a\|_{L^1(\mu)}\le C.
\end{aligned}
\]
This argument also applies to $F_j^*$ for every $1\le j\le n$.

It remains to estimate $R_1a$ on $\{x: |x-c|\ge3\}$.
Since $a$ has cancellation against $1,x_2,\ldots,x_n$,
\[
\begin{aligned}
 R_1a(x)=e^{-2x_1}\int_B
 \bigg[F_1(x-y)-F_1(x-c)
 +\sum_{k=2}^n(y_k-c_k)\partial_kF_1(x-c)\bigg]
 a(y)\,d\mu(y).
\end{aligned}
\]
The  Taylor expansion gives
\[
\begin{aligned}
 &F_1(x-y)-F_1(x-c)
   +\sum_{k=2}^n(y_k-c_k)\partial_kF_1(x-c)\\
 &=- (y_1-c_1)\int_0^1
       \partial_1F_1\bigl(x-y+t(y_1-c_1)e_1\bigr)\,dt\\
 &\quad+\sum_{k,\ell=2}^n(y_k-c_k)(y_\ell-c_\ell)
   \int_0^1(1-t)\,
       \partial_k\partial_\ell F_1
       \biggl(x-c-t\sum_{m=2}^n(y_m-c_m)e_m\biggr)\,dt.
\end{aligned}
\]
This, in combination with \eqref{eq:global-integrability},
gives
\[
\begin{aligned}
 \int_{|x-c|\ge3}|R_1a(x)|\,d\mu(x)
 &\le C\int_B\bigl(|y_1-c_1|+|y-c|^2\bigr)|a(y)|\,d\mu(y)\\
 &\le C(r+r^2)\|a\|_{L^1(\mu)}\le C.
\end{aligned}
\]

Hence $\|R_ja\|_{L^1(\mu)}+\|R_j^*a\|_{L^1(\mu)}\le C$. 

This completes our proof.
\end{proof}

As a direct consequence of Theorems \ref{thm:H1-endpoint}, \ref{thm:duality} and \ref{thm:BMO-interpolation}, we have:
\begin{theorem}\label{thm:BMO-endpoint}
For every $1\le j\le n$,
\begin{equation*}
 \|R_jg\|_{{\BMO(\mu)}}+\|R_j^*g\|_{{\BMO(\mu)}}
 \le C\|g\|_{L^\infty(\mu)}.
\end{equation*}
 Moreover, for $1<p<\infty$,
\[
 \|R_jf\|_{L^p(\mu)}+\|R_j^*f\|_{L^p(\mu)}
 \le C_p\|f\|_{L^p(\mu)}.
\]
\end{theorem}

{The following theorem shows that $H^1(\mu)$ cannot be characterized} by the Riesz transforms.
{
	
\begin{theorem}\label{thm:failure-Riesz-characterization}
There exists $f\in L^1(\mu)\cap L^2(\mu)$ such that
\[
R_jf,\ R_j^*f\in L^1(\mu),\qquad 1\le j\le n,
\]
but $f\notin {H^1(\mu)}$. In particular, ${H^1(\mu)}$ cannot be
characterized by the 
Riesz transforms.
\end{theorem}

\begin{proof}
Choose $\varphi\in C_c^\infty(\R^n)$ such that
$\int\varphi\,d\mu\ne0$, and set $f=\cL^{1/2}\varphi$.
The semigroup $e^{-t\cL}$ is contractive on $L^1(\mu)$. Since
\[
\cL^{1/2}\varphi
=\frac1{2\sqrt\pi}\int_0^\infty
(\varphi-e^{-t\cL}\varphi)\,\frac{dt}{t^{3/2}},
\]
we have
\[
\begin{aligned}
\|f\|_{L^1(\mu)}
&\le C\int_0^1 t^{-1/2}\|\cL\varphi\|_{L^1(\mu)}\,dt
+C\int_1^\infty t^{-3/2}\|\varphi\|_{L^1(\mu)}\,dt\\
&\le C\bigl(\|\cL\varphi\|_{L^1(\mu)}
+\|\varphi\|_{L^1(\mu)}\bigr)<\infty.
\end{aligned}
\]
Moreover, $f\in L^2(\mu)$, and since $\cL$ commutes with
$\partial_j$ and $\partial_j^*$, we obtain
\[
R_jf=\partial_j\varphi,\qquad
R_j^*f=\partial_j^*\varphi
=(-\partial_j-2\delta_{1j})\varphi.
\]
Thus $R_jf,R_j^*f\in L^1(\mu)$ for $1\le j\le n$.

For $h\in L^1(\mu)$, define
\[
\widehat h_\mu(\xi)=\int_{\R^n}h(x)e^{-ix\cdot\xi}\,d\mu(x).
\]
If $a$ is an $(1,\infty;\cP)_1$ atom supported in $B(c,r)$, $r\le1$,
then $\int a\,d\mu=0$ and $\|a\|_{L^1(\mu)}\le1$. Hence
\[
\begin{aligned}
|\widehat a_\mu(te_1)|
&=\left|\int_B\bigl(e^{-itx_1}-e^{-itc_1}\bigr)a(x)\,d\mu(x)\right|\\
&\le |t|r\|a\|_{L^1(\mu)}\le |t|.
\end{aligned}
\]
It follows that  for every $h\in {H^1(\mu)}$, we have
\begin{equation}\label{eq-f hardy}
|\widehat h_\mu(te_1)|\le |t|\|h\|_{{H^1(\mu)}}.
\end{equation}
On the other hand, using
$e^{2x_1}\cL e^{-2x_1}=-\Delta+2\partial_1$, we obtain
\[
\widehat f_\mu(\xi)
=\bigl(|\xi|^2+2i\xi_1\bigr)^{1/2}
\widehat\varphi_\mu(\xi).
\]
Indeed, this follows from the preceding semigroup formula and
\[
\widehat{e^{-t\cL}\varphi}_\mu(\xi)
=e^{-t(|\xi|^2+2i\xi_1)}\widehat\varphi_\mu(\xi).
\]
Since $\widehat\varphi_\mu(0)=\int\varphi\,d\mu\ne0$,
\[
|\widehat f_\mu(te_1)|
=|t^2+2it|^{1/2}|\widehat\varphi_\mu(te_1)|
\simeq |t|^{1/2},\qquad t\to0.
\]
This contradicts \eqref{eq-f hardy}.

Thus $f\notin {H^1(\mu)}$.

This completes our proof.
\end{proof}
}

\begin{theorem}\label{thm:sharp-moment}
Let $\cQ$ be any proper subset of
$\{1,x_2,\ldots,x_n\}$. There is a function $a$ supported in a ball $B$ of radius less than $1$, 
with $\|a\|_\infty\le\mu(B)^{-1}$ and
\[
 \int a(x)p(x)\,d\mu(x)=0, \ \ \ \ \qquad p\in\cQ,
\]
such that $R_1a\notin L^1(\mu)$. Consequently, if the atoms defining
${H^1(\mu)}$ are required to annihilate only $\cQ$, the Riesz transform
$R_1$ is not bounded from that Hardy space to $L^1(\mu)$.
\end{theorem}

\begin{proof}

Choose $0\le\varphi\in C_c^\infty(B(0,1/4))$,
positive on $B(0,1/8)$.  
Since $\varphi$ is positive on a ball, the form
\[
\langle p,q\rangle_\varphi
=\int\varphi(y)p(y)q(y)\,dy
\]
defines an inner product on the real span of
$\{1,x_2,\ldots,x_n\}$. Indeed,
$\int\varphi(y)p(y)^2\,dy>0$ for every nonzero real
polynomial $p$ in this span. As $\cQ$ is a proper subset,
there is a nonzero real polynomial $\psi\in\cP$ such that
\[
\int\varphi(y)\psi(y)p(y)\,dy=0,
\qquad p\in\cQ.
\]

Moreover, since $\psi\ne0$, we have
\[
\int\varphi(y)\psi(y)^2\,dy>0.
\]
Thus at least one of the moments
$\displaystyle \int\varphi(y)\psi(y)p(y)\,dy$,
$p\in\{1,y_2,\ldots,y_n\}$, is nonzero.

For sufficiently small $c>0$, set
$a(y)=ce^{-2y_1}\varphi(y)\psi(y)$. Then, we have
\begin{enumerate}
	\item[(i)] $\supp a\subset B(0,1/4)$;
	\item[(ii)] $\|a\|_{L^\infty(\mu)}\lesssim 1$;
	\item[(iii)] $\displaystyle \int a(x)p(x)d\mu(x)=0$ for all $p\in \mathcal Q$.
\end{enumerate} 
{Hence, $a$ is a multiple of an $(1,\infty;\mathcal Q)_{1/4}$-atom.} Moreover,
\[
 (M_0,M_2,\ldots,M_n)
 =\left(\int a\,d\mu,\int y_2a\,d\mu,\ldots,
                         \int y_na\,d\mu\right)\ne0.
\]
 
Recall that $G(z)=g(|z|)$, where
\[
 g(r)=d_n\int_0^\infty t^{-(n+1)/2}e^{-t-r^2/(4t)}\,dt,
 \qquad d_n=\frac1{\sqrt\pi(4\pi)^{n/2}}.
\]
{Set $\nu=(n-1)/2$. Then}
\[
g(r)=A_n r^{-\nu}K_\nu(r),
\qquad {A_n=2^{\nu+1}d_n,}
\]
{where $K_\nu$} is the modified
Bessel function  of the second kind  defined by
\[
K_\nu(r)=\frac12\left(\frac r2\right)^\nu
\int_0^\infty t^{-\nu-1}
e^{-t-r^2/(4t)}\,dt,\qquad r>0.
\]

Using the identities
\[
 \frac{d}{dr}\bigl(r^{-\nu}K_\nu(r)\bigr)
 =-r^{-\nu}K_{\nu+1}(r),\qquad
 K_\alpha'(r)=-K_{\alpha+1}(r)+\frac\alpha rK_\alpha(r),
\]
we obtain
\[
 \begin{aligned}
 g'(r)&=-A_n r^{-\nu}K_{\nu+1}(r),\\
 g''(r)&=A_n r^{-\nu}
   \left(K_{\nu+2}(r)-\frac1rK_{\nu+1}(r)\right).
 \end{aligned}
\]
For each fixed $\alpha>0$, the standard asymptotic property
\[
 K_\alpha(r)=\sqrt{\frac\pi{2r}}e^{-r}
 \bigl(1+O(r^{-1})\bigr),\qquad r\to\infty,
\]
 implies that   as $r\to\infty$,
\begin{equation}\label{eq:G-radial-asympt}
\begin{aligned}
 g(r)&=c_ne^{-r}r^{-n/2}(1+O(r^{-1})),\\
 g'(r)&=-c_ne^{-r}r^{-n/2}(1+O(r^{-1})),\\
 g''(r)&=c_ne^{-r}r^{-n/2}(1+O(r^{-1})),
\end{aligned}
\qquad c_n=(2\pi)^{-n/2}.
\end{equation}
 
Recall that
\[
 F_1(z)=e^{z_1}\left(\frac{z_1}{|z|}g'(|z|)-g(|z|)\right).
\]
For $n\ge2$, write $z=(z_1,z_\perp)$, where
$z_\perp=(z_2,\ldots,z_n)$, and define
$\nabla_\perp=(\partial_2,\ldots,\partial_n)$.
Hence, we have, with $r=|z|$,
\[
 \nabla_\perp F_1(z)
 =e^{z_1}\frac{z_\perp}{r}
 \left(\frac{z_1}{r}g''(r)
 -\frac{z_1}{r^2}g'(r)-g'(r)\right).
\]
Now set $z=(s,\sqrt{s}\,w)$, with $w\in\R^{n-1}$ bounded.
As $s\to\infty$, we have
\[
 \begin{aligned}
 r&=s+\frac{|w|^2}{2}+O(s^{-1}),
 &\frac{s}{r}&=1+O(s^{-1}),\\
 e^{s-r}&=e^{-|w|^2/2}(1+O(s^{-1})),
 &r^{-n/2}&=s^{-n/2}(1+O(s^{-1})).
 \end{aligned}
\]
It follows from \eqref{eq:G-radial-asympt} that
\begin{equation}\label{eq-F1 gradient F1}
\begin{aligned}
 F_1(s,\sqrt{s}\,w)
 &=-2c_ne^{-|w|^2/2}s^{-n/2}(1+O(s^{-1})),\\
 \nabla_\perp F_1(s,\sqrt{s}\,w)
 &=\bigl(2c_nw+O(s^{-1})\bigr)
 e^{-|w|^2/2}s^{-(n+1)/2}.
\end{aligned}
\end{equation}
Both estimates are uniform on bounded sets containing $w$. The second
formula is omitted when $n=1$.

By \eqref{eq:Fj}, for $|x|>3$ we have
\[
 R_1a(x)=\int K_1(x,y)a(y)\,d\mu(y)
 =e^{-2x_1}\int F_1(x-y)a(y)\,d\mu(y).
\]
Since $d\mu(x)=e^{2x_1}\,dx$, it follows that
\begin{equation}\label{eq:R1-sharp-L1}
 \int_{|x|>3}|R_1a(x)|\,d\mu(x)
 =\int_{|x|>3}\left|\int F_1(x-y)a(y)\,d\mu(y)\right|dx.
\end{equation}
We show that the integral on the right is infinite.

Set $M_\perp=(M_2,\ldots,M_n)$ and
\[
 H(x)=\int F_1(x-y)a(y)\,d\mu(y).
\]
We first prove that, for $|x|>3$,
\begin{equation}\label{eq:R1-sharp-expansion}
 H(x)=M_0F_1(x)-M_\perp\cdot\nabla_\perp F_1(x)+E(x),
 \qquad E\in L^1(\{|x|>3\},dx).
\end{equation}
Writing $y=y_1e_1+y_\perp$, we have
\[
 F_1(x-y)-F_1(x-y_\perp)
 =-y_1\int_0^1\partial_1F_1(x-y_\perp-ty_1e_1)\,dt
\]
and
\[
\begin{aligned}
 F_1(x-y_\perp)
 &=F_1(x)-\sum_{k=2}^ny_k\partial_kF_1(x) +\sum_{k,\ell=2}^ny_ky_\ell
 \int_0^1(1-t)\partial_k\partial_\ell F_1(x-ty_\perp)\,dt.
\end{aligned}
\]
{Combining these identities and integrating}  with respect to $a(y)\,d\mu(y)$, we obtain \eqref{eq:R1-sharp-expansion}, where
\[
\begin{aligned}
	E(x)=&-\int a(y)y_1\int_0^1\partial_1F_1(x-y_\perp-ty_1e_1)\,dt\, d\mu(y)\\
	& \ \ +\sum_{k,\ell=2}^n\int a(y)y_ky_\ell
	\int_0^1(1-t)\partial_k\partial_\ell F_1(x-ty_\perp)\,dt\,d\mu(y).
\end{aligned}
\]

 Since $|x-ty_\perp|, |x-y_\perp-ty_1e_1|>1$ whenever $|y|<1/2$, $|x|\ge3$ and $t\in (0,1)$, using Lemma~\ref{lem:kernel-derivatives} and
Fubini's theorem,
\[
 \int_{|x|>3}|E(x)|\,dx
 \le C\int\bigl(|y_1|+|y_\perp|^2\bigr)|a(y)|\,d\mu(y)<\infty.
\]

Setting $P(x)=M_0F_1(x)-M_\perp\cdot\nabla_\perp F_1(x)$, if  $n\ge2$, then for $x=(s,\sqrt{s}\,w)$,
$|w|\le1$, by \eqref{eq-F1 gradient F1}, we have
\[
\begin{aligned}
 P(s,\sqrt{s}\,w)
 &=-2c_ne^{-|w|^2/2}M_0s^{-n/2}
   +O\bigl(|M_0|s^{-n/2-1}\bigr)\\
 &\quad-2c_ne^{-|w|^2/2}(M_\perp\cdot w)s^{-(n+1)/2}
   +O\bigl(|M_\perp|s^{-(n+3)/2}\bigr),
\end{aligned}
\]
uniformly in $w$.

If $M_0\ne0$,  $|P(s,\sqrt{s}\,w)|\ge c|M_0|s^{-n/2}$ for $|w|\le1$.
Using $dx=s^{(n-1)/2}\,ds\,dw$, we obtain
\[
 \int_{|x|>3}|P(x)|\,dx
 \ge c|M_0|\int_{s_0}^\infty s^{-1/2}\,ds=\infty.
\]
If $M_0=0$, then $M_\perp\ne0$. Choose a measurable set
$A\subset\{|w|\le1\}$ of positive measure such that
$|M_\perp\cdot w|\ge c|M_\perp|$ on $A$.
For sufficiently large $s$,  $|P(s,\sqrt{s}\,w)|\ge c|M_\perp|s^{-(n+1)/2}$ on $A$. Hence
\[
 \int_{|x|>3}|P(x)|\,dx
 \ge c|M_\perp|\int_{s_0}^\infty s^{-1}\,ds=\infty.
\]
If $n=1$, then we have $M_0\ne0$ and $P=M_0F_1$. Since
$F_1(s)=-2c_1s^{-1/2}(1+O(s^{-1}))$ as $s\to\infty$,
we also have $P\notin L^1(\{|x|>3\},dx)$.

Thus $P\notin L^1(\{|x|>3\},dx)$ in every case. Since
$H=P+E$ and $E\in L^1(\{|x|>3\},dx)$, it follows that
$H\notin L^1(\{|x|>3\},dx)$. By \eqref{eq:R1-sharp-L1},
$R_1a\notin L^1(\mu)$. 

This completes our proof.

\end{proof}

\section{Commutators of local BMO functions and Riesz transforms}\label{sec:commutators}

For a suitable function $b$, the commutator $[b,T]$ {of $b$ and a linear operator $T$} {is defined by} $[b,T]=bT-T(b\,\cdot)$.  In this section, we show that the BMO space defined in Definition \ref{def:BMO} {does not characterize the boundedness of the commutators of the Riesz transforms.}


\subsection{A counterexample for the Hardy-space dual BMO}

\begin{proposition}\label{prop:comm-counterexample}
	Suppose $n\ge2$ and let $b(x)=x_2^2$. Then $b\in{\BMO(\mu)}$ as defined in
	Section~\ref{sec:spaces}, but $[b,R_2]$ is unbounded on $L^2(\mu)$.
\end{proposition}

\begin{proof}
	For $B=B(c,r)$, $r\le1$, choose
	$P_B(x)=2c_2x_2-c_2^2\in\cP$. Then
	\[
	|b(x)-P_B(x)|=|x_2-c_2|^2\le r^2\quad \ \ \text{for all $x\in B$}.
	\]
	It follows that 
	\[
	\frac{1}{\mu(B)}\int_B |b(x)-P_B(x)| d\mu(x) \lesssim 1
	\]
	for all $B=B(c,r)$ with $r\le 1$.

	This, together with ~\eqref{eq:best-approx}, implies that
	$\|b\|_{{\BMO(\mu)}}\le C$.
	
	Use the unitary operator $Uf=e^{x_1}f$ from
	\eqref{eq:ground-state}. Since multiplication by $x_2$ commutes with $U$,
	\[
	U[b,R_2]U^{-1}=[x_2^2,T],
	\qquad T=\partial_2(-\Delta+1)^{-1/2}.
	\]
	The multiplier of $T$ is
	$m(\xi)=i\xi_2(1+|\xi|^2)^{-1/2}$. Hence
	$[x_2,T]$ is a   Fourier multiplier, with symbol
	\[
	i\partial_{\xi_2}m(\xi)
	=-\frac{1+\sum_{k\ne2}\xi_k^2}{(1+|\xi|^2)^{3/2}}.
	\]
	Let $\tau_a f(x)=f(x-ae_2)$. Translation invariance of $T$ gives
	\begin{equation}\label{eq:translate-commutator}
		\tau_{-a}[x_2^2,T]\tau_a
		=[x_2^2,T]+2a[x_2,T].
	\end{equation}
Choose $f\in{ C_c^\infty(\R^n)}$ with
	$[x_2,T]f\ne0$. The reverse triangle inequality in
	\eqref{eq:translate-commutator} yields
	\[
	\|[x_2^2,T]\tau_a f\|_{L^2(dx)}
	\ge2|a|\|[x_2,T]f\|_{L^2(dx)}-\|[x_2^2,T]f\|_{L^2(dx)}.
	\]
	The right side tends to infinity, whereas
	$\|\tau_a f\|_{L^2(dx)}=\|f\|_{L^2(dx)}$. This proves the required unboundedness.
	
	This completes our proof.
\end{proof}


\subsection{Characterization by commutators}

\begin{definition}\label{def:bmo-loc}
	The space $\BMO_{\mathrm{loc}}(\mu)$ is the set of all  locally integrable functions $b$ satisfying
	\begin{equation*}
		\|b\|_{\BMO_{\mathrm{loc}}(\mu)}
		=\sup_{B:\,r_B\le1}\frac1{\mu(B)}
		\int_B|b-b_{B,\mu}|\,d\mu,
		\qquad
		b_{B,\mu}=\frac1{\mu(B)}\int_B b\,d\mu.
	\end{equation*}
	{We require $\|b\|_{\BMO_{\mathrm{loc}}(\mu)}<\infty$.}
	We identify two symbols differing by a constant.
\end{definition}

{ 
\begin{lemma}\label{lem:local-BMO-scales}
For every fixed $a>0$, the seminorm
\[
\sup_{B:\,r_B\le a}\frac1{\mu(B)}
\int_B|b-b_{B,\mu}|\,d\mu
\]
is equivalent to $\|b\|_{\BMO_{\mathrm{loc}}(\mu)}$. 
Moreover, for all $1\le q<\infty$,
\begin{equation}\label{eq:loc-JN-Lq}
	\sup_{B:\,r_B\le a}
	\left(\frac1{|B|}\int_B|b-b_{B,dx}|^q\,dx\right)^{1/q}
	{\simeq_{a,q}\|b\|_{\BMO_{\mathrm{loc}}(\mu)}}.
\end{equation}
\end{lemma}
}

\begin{proof}
For $a>0$, let
\[
\beta_a=\sup_{B:\,r_B\le a}\frac1{\mu(B)}
\int_B|b-b_{B,\mu}|\,d\mu.
\]
Fix $0<a<A$. Clearly, $\beta_a\le\beta_A$. Let $B=B(c,r)$ with
$a<r\le A$. Then {we can cover $B$ by a family of balls} $\{B_i: i \in \mathcal I\}$ satisfying the following properties:
\begin{enumerate}
	\item[(i)] {$r_{B_i}=a$ and $B_i\subset 2B$ for all $i\in \mathcal I$;}
	\item[(ii)]  {the overlap graph admits a spanning tree $T$ such that $B_i\cap B_j$ contains a ball of radius $a/10$ for every edge $(i,j)$ of $T$;}
	\item[(iii)] {$\sharp \mathcal I \lesssim (A/a)^n$;}
	\item[(iv)] $\sum_{i\in \mathcal I}\mathbf{1}_{B_i}\lesssim 1$.
\end{enumerate}

For each  pair $B_i,B_j$ with $(i,j)\in T$, choose a ball
$D_{ij}\subset B_i\cap B_j$ of radius ${a/10}$. By local doubling property 
\eqref{eq-local doubling},
\[
\begin{aligned}
|b_{B_i,\mu}-b_{B_j,\mu}|
&\le \frac1{\mu(D_{ij})}\int_{D_{ij}}
\bigl(|b-b_{B_i,\mu}|+|b-b_{B_j,\mu}|\bigr)\,d\mu\\
&\lesssim_a\beta_a.
\end{aligned}
\]
At this stage, arguing similarly to the proof of Proposition \ref{prop:BMO-scale}, we obtain $\beta_A\lesssim_{a,A}\beta_a$, and hence
{$\beta_a\simeq_{a,A}\beta_A$ for every fixed $A>a>0$.}

Finally, if $B=B(c,r)$ with $r\le a$, then
\[
e^{-2a}e^{2c_1}\le e^{2x_1}\le e^{2a}e^{2c_1},
\qquad x\in B.
\]
Consequently,
\[
\inf_{m\in\mathbb C}\frac1{\mu(B)}\int_B|b-m|\,d\mu
\simeq_a
\inf_{m\in\mathbb C}\frac1{|B|}\int_B|b-m|\,dx.
\]
Since
\[
\inf_{m\in\mathbb C}\frac1{\mu(B)}\int_B|b-m|\,d\mu
\le\frac1{\mu(B)}\int_B|b-b_{B,\mu}|\,d\mu
\le 2\inf_{m\in\mathbb C}\frac1{\mu(B)}\int_B|b-m|\,d\mu,
\]
we obtain, for any fixed $a>0$,
\begin{equation*}
	\|b\|_{\BMO_{\mathrm{loc}}(\mu)}
	\simeq_a
	\sup_{B:\,r_B\le a}\inf_{c\in\mathbb C}
	\frac1{|B|}\int_B|b-c|\,dx.
\end{equation*}
{This proves \eqref{eq:loc-JN-Lq} for $q=1$.} 

The equivalence \eqref{eq:loc-JN-Lq} for $q>1$  follows from the standard dyadic proof of the John--Nirenberg
inequality on bounded cubes, followed by the covering argument above. Hence, we omit the details.

This completes our proof.

\end{proof}

The following theorem proves that the BMO functions defined in {Definition~\ref{def:bmo-loc}} are the right space for the boundedness of the commutators of the Riesz transforms. 
\begin{theorem}\label{thm:comm-characterization}
	Let $1<p<\infty$, $1\le j\le n$, and $b\in L^1_{\mathrm{loc}}(\R^n)$.
	Then the following are equivalent:
	\begin{enumerate}[label=(\roman*)]
		\item[\rm (i)] $b\in\BMO_{\mathrm{loc}}(\mu)$;
		\item[\rm (ii)] $[b,R_j]$ is bounded on $L^p(\mu)$;
		\item[\rm (iii)] $[b,R_j^*]$ is bounded on $L^p(\mu)$.
	\end{enumerate}
	Moreover, there exists $C_{n,p}>1$ such that 
	\begin{equation}\label{eq:comm-equivalence}
		C_{n,p}^{-1}\|b\|_{\BMO_{\mathrm{loc}}(\mu)}
		\le \|[b,R_j]\|_{L^p(\mu)\to L^p(\mu)}
		\le C_{n,p}\|b\|_{\BMO_{\mathrm{loc}}(\mu)},
	\end{equation}
	and the same two-sided estimate holds with $R_j^*$ in place of $R_j$.
	
\end{theorem}

Note that for $1<p<\infty$,
\[
[b,R_j^*]=-[\overline b,R_j]^*,\qquad
\|[b,R_j^*]\|_{L^p(\mu)\to L^p(\mu)}
=\|[\overline b,R_j]\|_{L^{p'}(\mu)\to L^{p'}(\mu)}.
\]
{Thus, (ii) and (iii) are equivalent.}  It suffices to prove (i) $\Rightarrow$ (ii)  and  (ii) $\Rightarrow$ (i), and \eqref{eq:comm-equivalence} for $R_j$.

We first prove sufficiency using the commutator theorem of
Coifman, Rochberg and Weiss~\cite{CRW}. We then prove necessity
by testing on two separated balls.

\begin{lemma}\label{lem:local-bmo-extension}
	Fix the unit cubes $Q_k=k+[0,1)^n$, $k\in\mathbb Z^n$. For
	$b\in\BMO_{\mathrm{loc}}(\mu)$ there exist functions
	$b_k\in\BMO(\R^n,dx)$ such that
	\begin{equation*}
		b_k=b-c_k\quad\text{on }Q_k+[-2,2]^n,
		\qquad
		\|b_k\|_{\BMO(dx)}\le C\|b\|_{\BMO_{\mathrm{loc}}(\mu)},
	\end{equation*}
	where $c_k$ is a constant and $C$ is independent of $k$.
\end{lemma}

\begin{proof}
	Put $\beta=\|b\|_{\BMO_{\mathrm{loc}}(\mu)}$.
	Let $Q_k^*=k+[-4,5]^n$, set
	$c_k=|Q_k^*|^{-1}\int_{Q_k^*}b\,dx$, and choose
	$\eta_k\in C_c^\infty(k+(-4,5)^n)$ equal to one on
	$k+[-2,3]^n$, with derivatives uniformly bounded in $k$.
	Set $b_k=\eta_k(b-c_k)$.
	
Then by \eqref{eq-local doubling} and \eqref{eq:loc-JN-Lq}, we have
	\begin{equation}\label{eq:extension-large-L1}
		\begin{aligned}
			\int_{Q_k^*}|b-c_k|\,dx&\simeq \frac{1}{|Q_k^*|}\int_{Q_k^*}|b-c_k|\,dx\\
			&\le 2\frac{1}{|Q_k^*|}\int_{Q_k^*}|b-b_{B(k,10\sqrt n),dx}|\,dx\\
			&\lesssim \frac{1}{|B(k,10\sqrt n)|}\int_{B(k,10\sqrt n)}|b-b_{B(k,10\sqrt n),dx}|\,dx\\
			&\lesssim \beta.
		\end{aligned}
	\end{equation}
	{Hence, if $B$ is a ball} and its radius is greater than or equal to $1$, then by \eqref{eq:extension-large-L1},
	\[
	\frac1{|B|}\int_B|b_k|\,dx\le C\beta,
	\]
	which implies
	\[
	\sup_{B: \ r_B\ge 1}\frac1{|B|}\int_B|b_k-(b_k)_{B,dx}|\,dx\le C\beta.
	\]
We now consider the estimate on small balls. Let $B=B(x_B,r_B)$ with $r_B<1$. If $B\cap\supp\eta_k=\emptyset$, then
$b_k=0$ on $B$. Otherwise, choose $m\ge1$ such that
$1\le 2^mr_B<2$, and set $B_j=B(x_B,2^jr_B)$, $0\le j\le m$.
By \eqref{eq:loc-JN-Lq},
\[
\begin{aligned}
|b_{B_j,dx}-b_{B_{j+1},dx}|
&\le \frac1{|B_j|}\int_{B_j}|b-b_{B_{j+1},dx}|\,dy\\
&\le 2^n\frac1{|B_{j+1}|}
\int_{B_{j+1}}|b-b_{B_{j+1},dx}|\,dy
\le C\beta,
\end{aligned}
\]
for $0\le j<m$.
{Since $B$ meets $\supp\eta_k$ and $r_B<1$},
$B_m\cup Q_k^*\subset D_k:=B(k,10\sqrt n)$. Hence,
using \eqref{eq:loc-JN-Lq} again,
\[
\begin{aligned}
|b_{B_m,dx}-c_k|
&\le |b_{B_m,dx}-b_{D_k,dx}|
     +|b_{D_k,dx}-b_{Q_k^*,dx}|\\
&\le \left(\frac{|D_k|}{|B_m|}
        +\frac{|D_k|}{|Q_k^*|}\right)
\frac1{|D_k|}\int_{D_k}|b-b_{D_k,dx}|\,dy
\le C\beta.
\end{aligned}
\]
Therefore,
\[
|b_{B,dx}-c_k|\le C(m+1)\beta
{\le C\beta\bigl(1+\log(1/r_B)\bigr)}.
\]
{Put $c=\eta_k(x_B)(b_{B,dx}-c_k)$.} Then, for $y\in B$,
\[
b_k(y)-c
=\eta_k(y)(b(y)-b_{B,dx})
 +\bigl(\eta_k(y)-\eta_k(x_B)\bigr)(b_{B,dx}-c_k).
\]
Since $\|\eta_k\|_\infty+\|\nabla\eta_k\|_\infty\le C$, we obtain
\[
\begin{aligned}
\inf_{c\in\mathbb C}\frac1{|B|}\int_B|b_k-c|\,dy
&\le C{\frac1{|B|}\int_B|b-b_{B,dx}|\,dy
	   +Cr_B|b_{B,dx}-c_k|}\\
&\le C\beta{\bigl[1+r_B(1+\log(1/r_B))\bigr]}
\le C\beta.
\end{aligned}
\]
It follows that
\[
\sup_{B: \ r_B<1}\frac1{|B|}\int_B|b_k-(b_k)_{B,dx}|\,dy\le C\beta.
\]

	This completes our proof.
\end{proof}

\begin{proof}[Proof of Theorem~\ref{thm:comm-characterization}: $(i)\Rightarrow(ii)$]

	Fix $1<p<\infty$, $1\le j\le n $, and let $V_pf(x)=e^{2x_1/p}f(x)$. Then $V_p$ is an
	isometry from $L^p(\mu)$ onto $L^p(dx)$. Using $\Phi(z)=e^{z_1}G(z)$ and \eqref{eq:Fj}, we find that
	$V_pR_jV_p^{-1}$ is a convolution operator with kernel
	\begin{equation}\label{eq:comm-conjugated-kernels}
		H(z)=e^{(2/p-1)z_1}
		(\partial_jG(z)-\delta_{1j}G(z)).
	\end{equation}
	Since $|2/p-1|<1$, arguing similarly to the  proof of Lemma~\ref{lem:test-Riesz-bounded}, we have
	\begin{equation}\label{eq:comm-kernel-local-global}
		\begin{aligned}
			|H(z)|+|z| |\nabla H(z)|&\le C_p|z|^{-n}, &&0<|z|\le2,\\
			|H(z)|+|\nabla H(z)|&\le C_p e^{-c_p|z|}, &&|z|\ge1.
		\end{aligned}
	\end{equation}
	
	Choose $\chi\in C_c^\infty(B(0,2))$ such that {$\chi=1$} on $B(0,1)$,
	and write $H=H_{\mathrm{loc}}+H_{\mathrm{glo}}$, where
	$H_{\mathrm{loc}}=\chi H$. The operator
	$T_{\mathrm{loc}}f=H_{\mathrm{loc}}*f$ is a standard
	Calder\'on--Zygmund operator: it satisfies the local size and smoothness
	bounds in \eqref{eq:comm-kernel-local-global}, and it is bounded on
	$L^2(dx)$. To prove the $L^2(dx)$ boundedness, we write
\[
 H_{\mathrm{loc}}=\partial_jG+K,
\]
where
\[
\begin{aligned}
 K(z)={}&\chi(z)(e^{(2/p-1)z_1}-1)\partial_jG(z)
       -\delta_{1j}\chi(z)e^{(2/p-1)z_1}G(z) -(1-\chi(z))\partial_jG(z).
\end{aligned}
\]
We show that $K\in L^1(dx)$. From the integral representation of $G$, we have, for $0<|z|\le2$,
\[
 |\nabla G(z)|\le C|z|^{-n},\qquad
 |G(z)|\le C
 \begin{cases}
 |z|^{1-n},&n\ge2,\\
 1+|\log|z||,&n=1.
 \end{cases}
\]
Since $|e^{(2/p-1)z_1}-1|\le C_p|z|$ for $|z|\le2$, it follows that
\[
\begin{aligned}
 \int_{|z|\le2}|\chi(z)|
 &\bigl(|e^{(2/p-1)z_1}-1|\,|\partial_jG(z)|+e^{(2/p-1)z_1}|G(z)|\bigr)\,dz\\
 &\le C_p\int_{|z|\le2}(|z|^{1-n}+|G(z)|)\,dz<\infty.
\end{aligned}
\]
Moreover, $1-\chi=0$ on $B(0,1)$. By \eqref{eq:Bessel-tail},
\[
 \int_{\R^n}|1-\chi(z)|\,|\partial_jG(z)|\,dz
 \le C\int_{|z|\ge1}e^{-|z|}|z|^{-n/2}\,dz<\infty.
\]
Hence $K\in L^1(dx)$. By the Fourier multiplier formula for $G$,
\[
 \widehat{\partial_jG}(\xi)
 =i\xi_j(1+|\xi|^2)^{-1/2}.
\]
Thus Plancherel's theorem and Young's inequality give
\[
 \|T_{\mathrm{loc}}f\|_{L^2(dx)}
 \le\|\partial_jG*f\|_{L^2(dx)}+\|K*f\|_{L^2(dx)}
 \le(1+\|K\|_{L^1(dx)})\|f\|_{L^2(dx)}.
\]
	
	Let $\beta=\|b\|_{\BMO_{\mathrm{loc}}(\mu)}$.
	For each $k\in\mathbb Z^n$ put
	$f_k=f\,\mathbf1_{Q_k+[-2,2]^n}$ and choose $b_k$ from
	Lemma~\ref{lem:local-bmo-extension}. If $x\in Q_k$, then
	$T_{\mathrm{loc}}f(x)=T_{\mathrm{loc}}f_k(x)$, and
	$b_k=b-c_k$ on $Q_k+[-2,2]^n$. Thus
	\[
	[b,T_{\mathrm{loc}}]f(x)
	=[b_k,T_{\mathrm{loc}}]f_k(x),\qquad x\in Q_k.
	\]
	Applying the Coifman--Rochberg--Weiss commutator estimate in ~\cite{CRW}, we have
	\begin{align}\label{eq:comm-local-estimate}
		\|[b,T_{\mathrm{loc}}]f\|_{L^p(dx)}^p
		&\le\sum_{k\in\mathbb Z^n}
		\|[b_k,T_{\mathrm{loc}}]f_k\|_{L^p(dx)}^p\notag\\
		&\le C_p\beta^p\sum_k\|f_k\|_{L^p(dx)}^p
		\le C_p\beta^p\|f\|_{L^p(dx)}^p.
	\end{align}
	
	We next estimate the global part. Let $m_k=|Q_k|^{-1}\int_{Q_k}b\,dx$. Local John--Nirenberg's inequality \eqref{eq:loc-JN-Lq} and the connectedness of the family $\{Q_k\}$ imply that  for every $1\le q<\infty$,
	\begin{equation}\label{eq:cube-bmo-estimates}
		\|b-m_k\|_{L^q(Q_k,dx)}\le C_q\beta,
		\qquad |m_k-m_\ell|\le C\beta(1+|k-\ell|).
	\end{equation}
	By \eqref{eq:comm-kernel-local-global},
	\[
	\sup_{x\in Q_k,\,y\in Q_\ell}
	|H_{\mathrm{glo}}(x-y)|
	\le C_p e^{-c_p|k-\ell|}.
	\]
Decompose
	\[
	b(x)-b(y)
	=(b(x)-m_k)+(m_k-m_\ell)+(m_\ell-b(y)).
	\]
	Using H\"older's inequality on the two unit cubes and
	\eqref{eq:cube-bmo-estimates} with exponents $p$ and $p'$, we obtain
	\begin{equation*}
		\left\|\int_{Q_\ell}
		H_{\mathrm{glo}}(\,\cdot-y)
		(b(\,\cdot)-b(y))f(y)\,dy\right\|_{L^p(Q_k,dx)}
		\le C_p\beta(1+|k-\ell|)e^{-c_p|k-\ell|}
		\|f\|_{L^p(Q_\ell,dx)}.
	\end{equation*}
	Since $(1+|k|)e^{-c_p|k|}\in\ell^1(\mathbb Z^n)$,
	discrete Young's inequality yields
	\begin{equation}\label{eq:comm-far-estimate}
		\|[b,T_{\mathrm{glo}}]f\|_{L^p(dx)}\le C_p\beta\|f\|_{L^p(dx)}.
	\end{equation}

Combining \eqref{eq:comm-local-estimate} and
\eqref{eq:comm-far-estimate}, we obtain
\[
\|[b,T_{\mathrm{loc}}]f+[b,T_{\mathrm{glo}}]f\|_{L^p(dx)}
\le C_p\beta\|f\|_{L^p(dx)}.
\]

Using \eqref{eq:comm-conjugated-kernels}, we obtain 
\[
\begin{aligned}
V_p[b,R_j]V_p^{-1}f
&=b(V_pR_jV_p^{-1}f)-V_pR_jV_p^{-1}(bf)\\
&=b(T_{\mathrm{loc}}+T_{\mathrm{glo}})f
 -(T_{\mathrm{loc}}+T_{\mathrm{glo}})(bf)\\
&=[b,T_{\mathrm{loc}}]f+[b,T_{\mathrm{glo}}]f.
\end{aligned}
\]
Consequently, for $u=V_p^{-1}f$,
\[
\begin{aligned}
\|[b,R_j]u\|_{L^p(\mu)}
&=\|V_p[b,R_j]V_p^{-1}f\|_{L^p(dx)}\\
&\le C_p\beta\|f\|_{L^p(dx)}
=C_p\beta\|u\|_{L^p(\mu)}.
\end{aligned}
\]

	This completes our proof.
\end{proof}

\paragraph{Proof of necessity.}

\begin{lemma}\label{lem:nondegenerate-cone}
	{For each $1\le j\le n$, there is $c>0$ such that}
	\begin{equation*}
		{-F_j(z)\ge c|z|^{-n}}
		\quad\text{if }0<|z|\le 1
		\text{ and } z_j\ge\tfrac12|z|.
	\end{equation*}
	{The constants are independent of $j$.}
\end{lemma}

\begin{proof}
	Recall that $\Phi(z)=e^{z_1}G(z)$, where
	\[
	G(z)=\frac{(4\pi)^{-n/2}}{\sqrt\pi}
	\int_0^\infty t^{-(n+1)/2}
	e^{-t-|z|^2/(4t)}\,dt.
	\]
	Direct differentiation gives
	\begin{equation*}
		-\partial_jG(z)
		=\frac{(4\pi)^{-n/2}}{2\sqrt\pi}\,z_j
		\int_0^\infty t^{-(n+3)/2}
		e^{-t-|z|^2/(4t)}\,dt.
	\end{equation*}
	If $|z|\le 1$ and $z_j\ge |z|/2$, then 
	\[
	\begin{aligned}
		-\partial_jG(z)\gtrsim |z|\int_{|z|^2}^{2|z|^2} t^{-(n+3)/2}
		 \,dt \gtrsim |z|^{-n}.
	\end{aligned}
	\]
	
{ 
Since $G(z)\ge0$ and $e^{z_1}\ge e^{-1}$ for $|z|\le1$,
\[
\begin{aligned}
-F_j(z)&=e^{z_1}\bigl(-\partial_jG(z)+\delta_{1j}G(z)\bigr)\\
&{\ge e^{z_1}(-\partial_jG(z))\gtrsim |z|^{-n}.}
\end{aligned}
\]
}

	This completes our proof.
\end{proof}

\begin{proof}[Proof of Theorem~\ref{thm:comm-characterization}:
	$(ii)\Rightarrow(i)$]
	Since $R_j$ has a real kernel, we may assume that $b$ is a  real-valued function.
	Let $B=B(c,r)$ with $0<r\le1/8$, and put
	\[
	\widetilde B=B(c-6re_j,r).
	\]
	For $x\in B$ and $y\in\widetilde B$, we have
	\[
	4r\le(x-y)_j\le |x-y|\le8r\le1,
	\]
	which implies
	\[
	\qquad (x-y)_j\ge\tfrac12|x-y|.
	\]
	By Lemma~\ref{lem:nondegenerate-cone} and \eqref{eq:Fj},
	\begin{equation*}
	-K_j(x,y)=-e^{-2x_1}F_j(x-y)
	\ge c e^{-2x_1}r^{-n},\qquad (x,y)\in B\times\widetilde B.
	\end{equation*}
	Moreover, by \eqref{eq:local-volume},
	\begin{equation}\label{eq:paired-ball-normalization}
	\mu(\widetilde B)\simeq\mu(B)\simeq e^{2c_1}r^n,
	\qquad e^{-2x_1}r^{-n}\mu(\widetilde B)\simeq1
	\quad(x\in B).
	\end{equation}

	Let $m$ be a median of $b$ on $\widetilde B$
	with respect to $\mu$. That is,
	\[
	\mu\{y\in\widetilde B:b(y)>m\}
	\le\frac12\mu(\widetilde B),
	\qquad
	\mu\{y\in\widetilde B:b(y)<m\}
	\le\frac12\mu(\widetilde B).
	\]
	Define
	\[
	F_-=\{y\in\widetilde B:b(y)\le m\},\qquad
	F_+=\{y\in\widetilde B:b(y)\ge m\}.
	\]
	Both sets have measure at least $\mu(\widetilde B)/2$.
	If $x\in B$ and $b(x)\ge m$, the sign of the kernel is constant and
	$b(x)-b(y)\ge b(x)-m$ for $y\in F_-$. Therefore, for $x\in B$ with $b(x)\ge m$,
	\begin{align}\label{eq:median-positive}
		\bigl|[b,R_j]\mathbf1_{F_-}(x)\bigr|
		&=\left|\int_{F_-}(b(x)-b(y))K_j(x,y)\,d\mu(y)\right|\notag\\
		&\ge c e^{-2x_1}r^{-n}\mu(F_-)(b(x)-m)
		\ge c(b(x)-m).
	\end{align}
	Similarly, for $x\in B$ with $b(x)\le m$,
	\begin{equation}\label{eq:median-negative}
		\bigl|[b,R_j]\mathbf1_{F_+}(x)\bigr|
		\ge c(m-b(x)).
	\end{equation}
	For each $x\in B$, we have
	\[
	[b,R_j]\mathbf1_F(x)
	=\int_F(b(x)-b(y))K_j(x,y)\,d\mu(y).
	\]

	Integrating \eqref{eq:median-positive}--\eqref{eq:median-negative},
	using the $L^p(\mu)$ boundedness of the commutator and
	\eqref{eq:paired-ball-normalization}, we conclude that
	\begin{align}\label{eq:median-Lp-control}
		\int_B|b-m|^p\,d\mu
		&\le C\left(\|[b,R_j]\mathbf1_{F_-}\|_{L^p(\mu)}^p
		+\|[b,R_j]\mathbf1_{F_+}\|_{L^p(\mu)}^p\right)\notag\\
		&\le C\|[b,R_j]\|_{L^p(\mu)\to L^p(\mu)}^p\mu(\widetilde B)\\
		&\le C\|[b,R_j]\|_{L^p(\mu)\to L^p(\mu)}^p\mu(B).
	\end{align}
	Using H\"older's inequality and
	$|b_{B,\mu}-m|\le\mu(B)^{-1}\int_B|b-m|\,d\mu$, we have
	\begin{equation*}
		\frac1{\mu(B)}\int_B|b-b_{B,\mu}|\,d\mu
		\le \frac2{\mu(B)}\int_B|b-m|\,d\mu\le C\|[b,R_j]\|_{L^p(\mu)\to L^p(\mu)}.
	\end{equation*}
	{for all balls $B=B(c,r)$ with $0<r\le1/8$. Therefore,}
	\[
	\|b\|_{{\BMO_{\mathrm{loc}}(\mu)}}\le C\|[b,R_j]\|_{L^p(\mu)\to L^p(\mu)}.
	\]
{The lower estimate in \eqref{eq:comm-equivalence} follows.}

	This completes our proof.
\end{proof}

\section{General Laplacians with drift}\label{sec:general-v}
We now apply our results to the Laplacian with an arbitrary constant
drift. Let $v\in\R^n\setminus\{0\}$ and define
\[
\cL_v=-\Delta-2v\cdot\nabla,
\qquad d\mu_v(x)=e^{2\langle v,x\rangle}\,dx.
\]
The operator is nonnegative and self-adjoint on $L^2(\mu_v)$.
Set $\alpha=|v|$ and choose an orthogonal matrix $O$ such that
$O(v/\alpha)=e_1$. Using the following change of variables 
\begin{equation*}
 z=\alpha O(x),
 \qquad x=\alpha^{-1}O^T(z),
\end{equation*}
we have
\begin{equation}\label{eq:L-transform-general}
 \langle v,x\rangle=z_1,\qquad
 d\mu_v(x)=\alpha^{-n}d\mu(z),\qquad
 \cL_v=\alpha^2\cL.
\end{equation}

Set
\begin{equation*}
 \cP_v=\{c+\ell\cdot x:\ c\in\mathbb C,\ \ell\cdot v=0\}.
\end{equation*}
Let $1<q\le\infty$ and $b>0$. An $(1,q;\cP_v)_b$-atom is a function $a$
supported in a ball $B=B(c,r)$, $r\le b$, such that
\[
 \|a\|_{L^q(\mu_v)}\le\mu_v(B)^{1/q-1},\qquad \text{and} \  
 \int_B a(x)p(x)\,d\mu_v(x)=0,\quad p\in\cP_v.
\]
 {We define $H^1_v(\mu_v)$ as the space of all functions $f\in L^1(\mu_v)$ which can be written as} 
\[
 f=\sum_{k=1}^\infty\lambda_k a_k\quad\text{in }L^1(\mu_v),
 \qquad \sum_{k=1}^\infty|\lambda_k|<\infty,
\]
where $a_k$ are $(1,q;\cP_v)_b$-atoms, with norm
\[
 \|f\|_{H^1_v(\mu_v)}=\inf\sum_{k=1}^\infty|\lambda_k|.
\]
Similarly to Section 2, we can prove that the space and its norm are independent of $q\in (1,\infty]$ and $b>0$.

For $B=B(c,r)$, choose an orthonormal basis
$w_1,\ldots,w_{n-1}$ of $v^\perp=\{x\in\R^n:\langle x,v\rangle=0\}$ and define
\[
 q_0(x)=1,\qquad q_k(x)=\frac{w_k\cdot(x-c)}r,
 \quad 1\le k\le n-1.
\]
These functions form a basis for $\cP_v$. For $f\in L^1(B,\mu_v)$,
let $\Pi_{B,v}f\in\cP_v$ be determined by
\begin{equation}\label{eq:projection-v}
 \int_B(f-\Pi_{B,v}f)\overline p\,d\mu_v=0,
 \qquad p\in\cP_v.
\end{equation}
Indeed, writing $\Pi_{B,v}f=\sum_{k=0}^{n-1}a_kq_k$, we have
\[
 G^v a=d^v,\quad
 G^v_{jk}=\frac1{\mu_v(B)}\int_B q_jq_k\,d\mu_v,\quad
 d^v_j=\frac1{\mu_v(B)}\int_B f q_j\,d\mu_v.
\]
Since
\[
 \sum_{j,k=0}^{n-1}\xi_jG^v_{jk}\overline{\xi_k}
 =\frac1{\mu_v(B)}\int_B\Big|\sum_{j=0}^{n-1}\xi_jq_j\Big|^2
 \,d\mu_v>0\qquad(\xi\ne0),
\]
$G^v$ is invertible. Thus $\Pi_{B,v}f$ is well defined.
We define
\[
 \|f\|_{\BMO_v(\mu_v)}
 =\sup_{r_B\le1}\frac1{\mu_v(B)}
 \int_B|f-\Pi_{B,v}f|\,d\mu_v.
\]
We use $\BMO_v(\mu_v)/\cP_v$, and define
\[
 \|b\|_{\BMO_{\mathrm{loc},v}(\mu_v)}
 =\sup_{r_B\le1}\frac1{\mu_v(B)}
 \int_B|b-b_{B,\mu_v}|\,d\mu_v,\qquad
 b_{B,\mu_v}=\frac1{\mu_v(B)}\int_Bb\,d\mu_v.
\]

\begin{theorem}\label{thm:general-v}
For every $v\ne0$ the following statements hold.
\begin{enumerate}[label=(\roman*)]
\item $(H^1_v(\mu_v))^*\simeq\BMO_v(\mu_v)/\cP_v$.
\item {The space $\BMO_v(\mu_v)$ satisfies the John--Nirenberg
inequality on balls of radius at most $1$.
}\item For $1\le j\le n$, the Riesz transforms
$R_{v,j}=\partial_j\cL_v^{-1/2}$ and $R_{v,j}^*$ satisfy
\[
 \|R_{v,j}f\|_{L^1(\mu_v)}+\|R_{v,j}^*f\|_{L^1(\mu_v)}
 \le C_v\|f\|_{H^1_v(\mu_v)},
\]
and
\[
 \|R_{v,j}g\|_{\BMO_v(\mu_v)}
 +\|R_{v,j}^*g\|_{\BMO_v(\mu_v)}
 \le C_v\|g\|_{L^\infty(\mu_v)}.
\]
\item The Hardy-space and BMO-space interpolation theorems hold
with $H^1_v(\mu_v)$, $\BMO_v(\mu_v)$ and $L^p(\mu_v)$ in place of
the corresponding normalized spaces.
\item Let $1<p<\infty$ and $1\le j\le n$. For
$b\in L^1_{\mathrm{loc}}(\R^n)$, the commutators
$[b,R_{v,j}]$ and $[b,R_{v,j}^*]$ are bounded on $L^p(\mu_v)$
if and only if $b\in\BMO_{\mathrm{loc},v}(\mu_v)$. Moreover,
\[
 \|[b,R_{v,j}]\|_{L^p(\mu_v)\to L^p(\mu_v)}
 \simeq_{n,p,v}\|b\|_{\BMO_{\mathrm{loc},v}(\mu_v)},
\]
with the same equivalence for $[b,R_{v,j}^*]$.
\end{enumerate}
\end{theorem}

\begin{proof}
For $1\le p<\infty$, define
\begin{equation*}
 (T_{v,p}f)(z)=\alpha^{-n/p}f(\alpha^{-1}O^T(z)).
\end{equation*}
By \eqref{eq:L-transform-general},
\[
 \|T_{v,p}f\|_{L^p(\mu)}^p
 =\alpha^{-n}\int_{\R^n}|f(\alpha^{-1}O^T(z))|^p\,d\mu(z)
 =\|f\|_{L^p(\mu_v)}^p.
\]
If $B_z=\alpha O(B)$, then
\[
 \mu(B_z)=\alpha^n\mu_v(B),\qquad
 \operatorname{rad}(B_z)=\alpha r
\]
for all balls $B=B(c,r)$.

In particular, for an $(1,q;\cP_v)_1$ atom $a$ supported in $B$,
\[
 \|T_{v,1}a\|_{L^q(\mu)}
 =\alpha^{n(1/q-1)}\|a\|_{L^q(\mu_v)}
 \le\mu(B_z)^{1/q-1},
\]
and, for $p\in\cP$,
\[
 \int_{B_z}(T_{v,1}a)(z)p(z)\,d\mu(z)
 =\int_B a(x)p(\alpha O(x))\,d\mu_v(x)=0.
\]
{The map $p(z)\mapsto p(\alpha O(x))$ takes $\cP$ onto $\cP_v$.}
Thus
\[
 \|f\|_{H^1_v(\mu_v)}
 =\|T_{v,1}f\|_{H^{1}_{\alpha}(\mu)},
\]
where $H^1_\alpha(\mu)$ is the atomic Hardy space defined by using $(1,q;\cP)_\alpha$  atoms {in Definition~\ref{def:H1}}.

For $b^\sharp(z)=b(\alpha^{-1}O^T(z))$, \eqref{eq:projection-v}
gives
\[
 (\Pi_{B,v}b)^\sharp=\Pi_{B_z}b^\sharp.
\]
Consequently,
\[
 \|b\|_{\BMO_v(\mu_v)}
 =\|b^\sharp\|_{\BMO_\alpha(\mu)},\qquad
 \|b\|_{\BMO_{\mathrm{loc},v}(\mu_v)}
 =\|b^\sharp\|_{\BMO_{\mathrm{loc},\alpha}(\mu)}.
\]
{Here $\BMO_\alpha(\mu)=\BMO_\alpha^1(\mu)$, and $\BMO_{\mathrm{loc},\alpha}(\mu)$ denotes the local BMO space defined by taking the supremum over balls with radius at most $\alpha$.}
{Statements (i) and (ii) follow from the corresponding results in Section~\ref{sec:spaces}.}

{By \eqref{eq:L-transform-general}, {for $u_j=O(e_j)$}} we have
\begin{equation*}
 T_{v,p}R_{v,j}T_{v,p}^{-1}
 =\sum_{k=1}^nO_{kj}R_k
 =\partial_{u_j}\cL^{-1/2},\qquad
 T_{v,2}R_{v,j}^*T_{v,2}^{-1}
 =\sum_{k=1}^nO_{kj}R_k^*.
\end{equation*}
{The corresponding identity for the adjoints holds on other $L^p$ spaces}
where the operators are bounded. Hence
\[
 \|R_{v,j}f\|_{L^1(\mu_v)}
 \le\sum_{k=1}^n|O_{kj}|\|R_kT_{v,1}f\|_{L^1(\mu)}
 \lesssim_{\alpha}\|f\|_{H^1_v(\mu_v)}.
\]
The adjoint bound follows in the same way. The BMO bounds follow
from the preceding projection identity and
Theorem~\ref{thm:BMO-endpoint}. This proves (iii).
Applying $T_{v,p}$ to the interpolation theorems of Section~2
proves (iv).

It remains to prove (v). For $u=O(e_j)$, let
$R_u=\partial_u\cL^{-1/2}=\sum_{k=1}^nu_kR_k$. Since
$T_{v,p}$ commutes with multiplication after the change of variables,
\[
 T_{v,p}[b,R_{v,j}]T_{v,p}^{-1}
 =[b^\sharp,R_u]
 =\sum_{k=1}^nu_k[b^\sharp,R_k].
\]
The sufficiency follows from Theorem~\ref{thm:comm-characterization}:
\[
 \|[b^\sharp,R_u]\|_{L^p(\mu)\to L^p(\mu)}
 \lesssim_{n,p}\|b^\sharp\|_{\BMO_{\mathrm{loc}}(\mu)}.
\]
For necessity, the function $F_u=\sum_{k=1}^nu_kF_k$ satisfies
\[
 F_u(z)=e^{z_1}(\partial_uG(z)-u_1G(z)),\qquad
 \partial_uG(z)=-\frac{d_n(u\cdot z)}2
 \int_0^\infty t^{-(n+3)/2}e^{-t-|z|^2/(4t)}\,dt.
\]
By the integral formula above, when $u\cdot z\ge|z|/2$,
\[
 -\partial_uG(z)\gtrsim|z|^{-n},\qquad
 G(z)=o(|z|^{-n}),\qquad |z|\to0.
\]
Thus, for some $\varepsilon>0$ independent of $|u|=1$,
\[
 -F_u(z)\gtrsim|z|^{-n},\qquad
 0<|z|\le\varepsilon,\quad u\cdot z\ge|z|/2.
\]
For $B=B(c,r)$ with $r\le\varepsilon/8$, take
$\widetilde B=B(c-6ru,r)$. If $x\in B$ and $y\in\widetilde B$, then
\[
 (x-y)\cdot u\ge4r,\qquad
 4r\le|x-y|\le8r\le\varepsilon.
\]
Thus the kernel of $R_u$ has one sign and satisfies
\[
 -K_u(x,y)\gtrsim e^{-2x_1}r^{-n},
 \qquad x\in B,\quad y\in\widetilde B.
\]
Arguing similarly to the proof of
Theorem~\ref{thm:comm-characterization}, we have 
\[
 \sup_{r_B\le\varepsilon/8}\frac1{\mu(B)}
 \int_B|b^\sharp-(b^\sharp)_{B,\mu}|\,d\mu
 \lesssim_{n,p}\|[b^\sharp,R_u]\|_{L^p(\mu)\to L^p(\mu)}.
\]
By similar arguments used in Lemma \ref{lem:local-BMO-scales}, we obtain
\[
\|b\|_{\BMO_{\mathrm{loc},v}(\mu_v)}\lesssim  \|[b,R_{v,j}]\|_{L^p(\mu_v)\to L^p(\mu_v)}.
\]
This completes our proof.
\end{proof}

\section*{Acknowledgements}
The first author was supported by the Australian Research Council
Discovery Project DP260101083. During the preparation of this manuscript, the authors
used ChatGPT (OpenAI) to assist with language editing,
LaTeX preparation, and checking mathematical arguments.
The authors reviewed and edited the content as needed
and take full responsibility for the content of the manuscript.

\end{document}